\documentclass{article}[10pt]
\usepackage{amsmath, amssymb, amsthm, epsf, graphicx, amscd}

\theoremstyle{plain}
\newtheorem{theorem}{Theorem}[section]  %% Added by Doug
\newtheorem{corollary}[theorem]{Corollary}
\newtheorem{lemma}[theorem]{Lemma}

\newtheorem{conjecture}[theorem]{Conjecture}
\theoremstyle{definition}

\theoremstyle{definition}             %%% Doug added
\newtheorem{remark}[theorem]{Remark}
\theoremstyle{definition}         %%% Doug added

\theoremstyle{definition}         %%% Doug added

\numberwithin{equation}{section}

\newcommand{\bb}[1]{\mathbb #1}
\def\f #1{ {\mathfrak #1} }

\def\bs #1{ {\boldsymbol #1} }
\newcommand{\ord}{\text{ord}}

\def\c #1{ {\mathcal #1} }

\def\Msk #1{ {\c M_{#1}^{(k)}} }  %% M_a symmetric tensor k-times
\def\MskDual #1{ {\c M_{#1}^{(k)*}} }  %% M_a^* symmetric tensor k-times
\def\Fsk #1{ {\bar \alpha_k( #1 )} }   %% k-th symm tensor of Frob
\newcommand{\Ldual}{{L^{\dag *}}}
\newcommand{\hatR}{{ \widehat{R} }}

\def\Fsk #1{ {\bar \alpha_k( #1 )} }   %% k-th symm tensor of Frob

\begin{document}

\title{$L$-functions of symmetric powers of cubic exponential sums}
\author{C. Douglas Haessig}
\date{\today}
\maketitle

%%% ----------------------------------------------------------------------
\begin{abstract}
For each positive integer $k$, we investigate the $L$-function
attached to the $k$-th symmetric power of the $F$-crystal associated
to the family of cubic exponential sums of $x^3 + \lambda x$ where
$\lambda$ runs over $\overline{\bb F}_p$. We explore its
rationality, field of definition, degree, trivial factors,
functional equation, and Newton polygon. The paper is essentially
self-contained, due to the remarkable and attractive nature of
Dwork's $p$-adic theory.

A novel feature of this paper is an extension of Dwork's effective
decomposition theory when $k < p$. This allows for explicit
computations in the associated $p$-adic cohomology. In particular,
the action of Frobenius on the (primitive) cohomology spaces may be
explicitly studied.
\end{abstract}


\setcounter{tocdepth}{1}
\tableofcontents

%%% ----------------------------------------------------------------------
\section{Introduction}

This paper represents a continuation in the study of
$L$-functions attached to symmetric powers of families of
exponential sums. Our central object of study will be the
family of cubic exponential sums of $x^3 + \lambda x$. Similar
studies have considered the Legendre family of elliptic curves
\cite{AdolphHecke} \cite{DwHecke}, and the $n$-variable
Kloosterman family \cite{FuWan} \cite{pR86}.

Various approaches have been used to study these $L$-functions.
Dwork \cite{DwHecke} used the existence of a Tate-Deligne mapping
(excellent lifting) of the elliptic family to study the symmetric
powers; this line of investigation was continued by Adolphson
\cite{AdolphHecke} in his thesis. Another related method is the
symmetric power of the associated $F$-crystal. This approach was
explored by Robba \cite{pR86} who used index theory to calculate the
degrees of such functions coming from the family of one-variable
Kloosterman sums. A third approach, developed by Fu and Wan
\cite{FuWan}, used $\ell$-adic methods of Deligne and Katz to study
Kloosterman sums in $n$-variables. Since we are interested in the
$p$-adic properties of the zeros and poles, and we believe excellent
lifting does not exist for our family when $p \equiv -1$ mod$(3)$,
we take an approach similar to that of Robba's.

The $L$-function, denoted $M_k(T)$, attached to the $k$-th symmetric
power of the cubic family is defined as follows. Fix a prime $p \geq
5$ and let $\zeta_p \in \bb C_p$ be a primitive $p$-th root of
unity. Let $\lambda$ be an element of $\overline{\bb F}_p$ and
denote by $deg(\lambda) := [\bb F_p(\lambda): \bb F_p]$ its degree.
Define $q := p^{deg(\lambda)}$. The $L$-function associated to the
cubic exponential sums
\[
S_m(\lambda) := \sum_{x \in \bb F_{q^m}} \zeta_p^{Tr_{\bb F_{q^m} /
\bb F_p}(x^3 + \lambda x)} \qquad m=1, 2, 3, \ldots
\]
is well-known to be a quadratic polynomial with coefficients in
$\bb Z[\zeta_p]$:
\[
L(\lambda, T) := \exp \left( \sum_{m \geq 1} S_m(\lambda)
\frac{T^m}{m} \right) = (1-\pi_1(\lambda) T)(1-\pi_2(\lambda)T).
\]
The $L$-function of the $k$-th symmetric power of this family is
defined by:
\[
M_k(T) := \prod_{\lambda \in |\bb A^1|} \prod_{i=0}^k (1 -
\pi_1(\lambda)^i \pi_2(\lambda)^{k-i} T^{deg(\lambda)})^{-1}
\]
where $|\bb A^1|$ denotes the set of Zariski closed points of
$\bb A^1$. The main theorem of this paper is the following.

\begin{theorem}\label{T: MainTheorem}
Let $k$ be a positive integer.
\begin{enumerate}
\item $M_k(T)$ is a rational function with coefficients in $\bb
Z[\zeta_p]$. Further, if $p \equiv 1 \> mod(3)$, then the
coefficients of $M_k(T)$ lie in the ring of integers of the unique
subfield $L$ of $\bb Q(\zeta_p)$ with $[L: \bb Q] = 3$. However, if
$p \equiv -1 \> mod(3)$ then $M_k(T)$ has integer coefficients.

\item For $k$ odd, $M_k(T)$ is a polynomial and satisfies
    the functional equation
\begin{equation}\label{E: IntroFE}
M_k(T) = c T^\delta M_k\left( p^{-(k+1)} T^{-1} \right)
\end{equation}
where $c$ is a nonzero constant and $\delta := deg \> M_k$.
Furthermore, if $k$ is odd and $k < p$, then $M_k(T)$ has
degree $(k+1)/2$, and writing
\[
M_k(T) = 1 + c_1 T + c_2 T + \cdots + c_{(k+1)/2}
T^{(k+1)/2},
\]
a quadratic lower bound for the Newton polygon is given by:
\[
ord_p(c_m) \geq \frac{(p-1)^2}{3p^2}(m^2 + (k+1)m)
\]
for $m = 0, 1, \ldots, (k+1)/2$.
\item For $k$ even, we may factorize $M_k(T)$ as follows:
\[
M_k(T) =
\begin{cases}
N_k(T)\widetilde{M}_k(T) & \text{for $k$ even and $k <
2p$} \\
\frac{N_k(T)\widetilde{M}_k(T)}{Q_k(T)} & \text{for $k$ even
and $k \geq 2p$}
\end{cases}
\]
where $\widetilde{M}_k(T)$ has degree $\leq k$ and
satisfies the functional equation (\ref{E: IntroFE}), and
\begin{alignat*}{2}
N_k(T) &= (1 - p^{k/2} T)^{m_k} & &\quad\text{if }p \equiv 1 \text{
mod}(12) \\
&= (1 - (-\bar g)^{k/2} T)^{m_k} & &\quad\text{if }p \equiv
5 \text{
mod}(12) \\
&= (1 - p^{k/2} T)^{n_k} (1 + p^{k/2} T)^{m_k - n_k} &
&\quad\text{if } 4 | k \text{ and } p \equiv 7 \text{ mod}(12) \\
&= (1 - p^{k/2} T)^{m_k - n_k} (1 + p^{k/2} T)^{n_k} &
&\quad\text{if } 4 \nmid k \text{ and } p \equiv 7 \text{ mod}(12) \\
&= (1 - \bar g^{k/2} T)^{n_k} (1 + \bar g^{k/2} T)^{m_k -
n_k} &
&\quad\text{if } 4 | k \text{ and } p \equiv 11 \text{ mod}(12) \\
&= (1 - \bar g^{k/2} T)^{m_k-n_k} (1 + \bar g^{k/2}
T)^{n_k} & &\quad\text{if }4 \nmid k \text{ and } p \equiv
11 \text{ mod}(12)
\end{alignat*}
where
\[
m_k :=
\begin{cases}
1 + \left\lfloor \frac{k}{2p} \right\rfloor & \text{if } 4|k \\
\left\lfloor \frac{k}{2p} \right\rfloor & \text{if } 4 \nmid k
\end{cases}
\quad \text{and} \quad n_k :=
\begin{cases}
1 + \left\lfloor \frac{k}{4p} \right\rfloor & \text{if } 4|k \\
\left\lfloor \frac{k}{4p} \right\rfloor & \text{if } 4 \nmid k
\end{cases}
\]
and where $\bar g$ is the complex conjugate of the Gauss
sum $g := g_2((p^2-1)/3)$ (see \S \ref{S: DeformTheory}).
Note, we conjecture $Q_k(T)$ always equals 1; see the
conjecture below.
\end{enumerate}
\end{theorem}

The lower bound for the Newton polygon is obtained using an
extension of Dwork's ``General Theory'' (see
\cite{DwZetaIHES}). This theory produces an explicit algorithm
to compute elements in cohomology. More specifically, given an
analytic element $\xi$ on the level of the cochain ($\c
M_a^{(k)}(b',b)$ in our notation) with known growth rate,
Theorem \ref{T: oddDecomp} describes the growth rate for the
reduction $\bar \xi$ of this element in the cohomolgy $H^1_k$.
Since we are able to easily understand the Frobenius action on
the cochain level in terms of growth rates, we can translate
this into growth rates in cohomology. The lower bound of the
Newton polygon follows (see Theorem \ref{T: OddNewtonPolygon}).

We expect the lower bound may be improved to $\frac{1}{3}(m^2 +
m + mk)$. This would be optimal since, for $k=3$ and $p=7$,
Vasily Golyshev has computed $ord_p(c_1) = 5/3$ which coincides
with this predicted lower bound. This expected lower bound
follows from the general philosophy that as $p \rightarrow
\infty$, the Newton polygon equals the Hodge polygon (see
\cite[Conjecture 1.9]{WanVariationNP}
\cite{ZhuAsympVariation}). (Note, the $\frac{k}{3}m$ term is
just an artifact leftover from the divisibility of the
Frobenius matrix by $p^{k/3}$.) We believe this lower bound
will also hold for $k$ even and $k < p$.

Looking for lower bounds of this type are inspired by a
reciprocity theorem of Wan's \cite{WanModular}. The
Gov\^{e}a-Mazur conjecture predicts a certain (essentially
linear) upper bound related to the arithmetic variation of
dimensions associated to classical modular forms of weight $k$
on $\Gamma_1(N) \cap \Gamma_0(p)$, with $(N,p)=1$. To prove a
quadratic upper bound for this conjecture, Wan
\cite{WanModular} proves a reciprocity theorem which cleverly
transforms a quadratic uniform lower bound for the Newton
polygon associated to an Atkin's operator on a space of
$p$-adic modular forms into a quadratic upper bound for the
Gov\^{e}a-Mazur conjecture.

At this time, we are only able to prove a lower bound for the Newton
polygon when $k$ is odd and $k < p$. The reason for the second
restriction comes from certain rational numbers in the proof of our
decomposition theorem (Theorem \ref{T: oddDecomp}). Their
denominators are $p$-adic units when $k < p$, but are often not
units when $k \geq p$.

We believe the obstruction to proving a similar bound for $k$ even
and $k < p$ lies in the underlying variety, or motive, associated to
the $k$-th symmetric power, being singular. Evidence for this is
given in Livn\'{e} \cite{LivneCubic}. In that paper, Livn\'{e}
relates the $k$-th moment of the cubic family $a x^3 + b x$ to $\bb
F_p$-rational points on a hypersurface $W_k$ in $\bb P^{k-2}$. When
$k$ is odd, $W_k$ is nonsingular, but when $k$ is even, $W_k$ has
ordinary double points. Livn\'{e} overcomes this difficulty by
embedding $W_k$ in a non-singular family. Interestingly, Dwork
proceeds in a similar manner when trying to generalize his effective
decomposition theory \cite[\S 3d]{DwZetaIHES}, which works well for
non-singular projective hypersurfaces, to singular ones
\cite{DwIII}. We are curious whether one can adapt these methods to
the current situation.

A first step in this direction may lie in a partial effective
decomposition theorem (Theorem \ref{T: evenDecomp}) proven when
$k$ is even and $k < p$. Even though the theorem does not give
complete information about the cohomology, it does produce
several non-trivial basis vectors of the cohomology and some
information concerning the lower bound of the Newton polygon.

Another topic is the denominator of $M_k(T)$. As we shall see
in \S \ref{S: cohomForm} and \S \ref{S: TrivialSubspaces},
$M_k(T)$ is a quotient of characteristic polynomials of
Frobenius matrices acting on cohomology spaces:
\begin{equation}\label{E: IntroCF}
M_k(T) = \frac{N_k(T) det(1 - \bar \beta_k T | PH_k^1)}{det(1-p \bar \beta_k
T | H_k^0)}.
\end{equation}
The polynomial $Q_k(T)$ is defined as $det(1-p \bar \beta_k T |
H_k^0)$. For every odd $k$, or for every even $k$ with $k <
2p$, $H_k^0$ has dimension zero. We suspect this is always
true:

\begin{conjecture}\label{C: NoOverConv}
$H_k^0 = 0$ for all positive integers $k$. Consequently,
$M_k(T)$ is a polynomial for all positive integers $k$.
\end{conjecture}

Our methods also allow us to study the $L$-function, denoted $M_k(d,
T)$, attached to the symmetric powers of the family $x^d + \lambda
x$, when $d$ is not divisible by $p$. In particular, we are able to
prove:

\begin{theorem}
$M_k(d,T)$ is a rational function with coefficients in $\bb
Z[\zeta_p]$. If $r := (d,p-1)$ then the coefficients of $M_k(d,T)$
lie in the ring of integers of the unique subfield $L$ of $\bb
Q(\zeta_p)$ with $[L: \bb Q] = r$. In particular, if $(d,p-1) = 1$
then $M_k(d,T)$ has integer coefficients.
\end{theorem}

We begin the paper by defining a relative cohomology theory (the
so-called $d$-Airy $F$-crystal) tailored specifically for the family
$x^d + a x$, where $a$ is the parameter. This will be a free module
$\c M_a(b',b)$ of rank $d$ over a power series ring $L(b')$. It will
carry an action of Frobenius $\alpha(a)$ and a connection
$\partial_a$. As we shall see in \S \ref{S: Fibres}, when we
specialize the parameter $a$ to a lifting in $\overline{\bb Q}_p$ of
an element $\bar z \in \overline{\bb F}_p^*$, the relative
cohomology reduces to a $\bb C_p$-vector space $M_z$, and the
$L$-function of $x^d+ \bar z x$ has the cohomological description as
the characteristic polynomial of the Frobenius on $M_z$:
\[
L(x^d+\bar z x, T) = det( 1 - \bar \alpha_{z,s} T | M_z)
\]
We then prove explicitly the functional equation for this
$L$-function via Dwork's theory. This means defining a relative dual
space $\c R'_{-\pi, a}(b',b)$ to $\c M_{\pi, a}(b',b)$, and an
isomorphism $\overline{\Theta}_{-\pi, a}$ between them which relates
the Frobenius $\bar \alpha_\pi(a)$ to its conjugate dual operator
$\bar \alpha_{-\pi}^*(a)$:
\[
\overline{\Theta}_{-\pi, a^p} = p^{-1} \bar \alpha_\pi(a) \circ
\overline{\Theta}_{-\pi, a} \circ \bar \alpha_{-\pi}^*(a).
\]
Lastly, we find the exact $p$-adic order of the entries in the
Frobenius matrix $\bar \alpha(a)$. This allows us to explicitly
describe the Newton polygon of $L(x^d+\bar z x, T)$; see \S
\ref{S: Fibres}.

Next, for $d=3$, the connection $\partial_a$ on the relative
cohomology produces a system of differential equations, the
Airy differential system, which has only an irregular singular
point at infinity. The (dual) Frobenius is an isomorphism on
the local solutions of this system and so, locally, is given by
an invertible constant matrix $M$. For local solutions near 0,
$M$ may be explicitly described in terms of Gauss sums. When
the solutions are near infinity, $M$ may be explicitly
described in terms of the square roots of $p$ or Gauss sums
$g$; see \S \ref{S: AiryAsymp} for details. The description of
$M$ depends on the congruence class of $p$ modulo 12 and is the
main reason for so many different forms of $N_k(T)$ in Theorem
\ref{T: MainTheorem}.

In \S \ref{S: cohomForm}, we present the general theory for a
cohomological formula for $M_k(T)$ following the work of Robba.
In short, the $k$-th symmetric power of the free module $\c
M_a(b',b)$ is another free module $\c M_a^{(k)}(b',b)$ which
carries a new action of Frobenius $\beta_k$ which is built from
the $k$-th symmetric power of the Frobenius $\bar \alpha(a)$,
and a new connection which is also denoted $\partial_a$. This
means we are able to take cohomology once again, creating the
finite dimensional $\bb C_p$-vector spaces $H_k^0$ and $H_k^1$.
The cohomological formula for $M_k(T)$ is equation (\ref{E:
IntroCF}) above.

In \S \ref{S: TrivialSubspaces}, we shall see that the dual space of
$H_k^1$ splits into three parts: a constant subspace $\bb
C_p^{k+1}$, a trivial subspace $\f T_k$, and a primitive part
$PH_k^{1*}$. For the first two, the action of Frobenius is
explicitly described. These descriptions lead to the polynomials
$P_k(T)$ and $N_k(T)$ described in Theorem \ref{T: MainTheorem}. The
action of the Frobenius on $PH_k^{1*}$ is more difficult to
understand, however, we are able to present the theory for the
functional equation of $\widetilde{M}_k(T)$. This is similar in
nature to that of $L(x^d+\bar z x, T)$, yet different since the
analogous operator to that of $\overline{\Theta}_{-\pi, a}$ has a
kernel which must be dealt with.

In \S \ref{S: CohomSymPow}, we present an effective
decomposition theory for the cohomology space $H_k^1$ when $k$
is odd and $k < p$. More precisely, an explicit procedure is
described which takes an element $\xi$ of $\c M_a^{(k)}(b)$ and
produces its reduction $\bar \xi$ in $H_k^1 = \c M_a^{(k)}(b) /
\partial_a \c M_a^{(k)}(b)$. As an application, we may compute a
non-trivial lower bound for the entries of the Frobenius matrix
$\bar \beta_k$ acting on $PH_k^1$. This produces the quadratic lower
bound for the Newton polygon given in Theorem \ref{T: MainTheorem};
see \S \ref{S: NPofM}.

\bigskip\noindent{\bf Acknowledgements.} I wish to thank Alan
Adolphson for many helpful comments, Steve Sperber for useful
conversations, Zhi-Wei Sun for his proof of the combinatorial
formula (\ref{E: Combo}), Vasily Golyshev for pointing out some
technical errors, and Daqing Wan for many insightful
conversations. Lastly, I would like to thank the referees for
suggesting many improvements and corrections.

%%% ----------------------------------------------------------------------
\section{Relative Dwork Theory}\label{S: RDT}

%%% ----------------------------------------------------------------------
\subsection{Relative Cohomology}\label{SS: RelCohom}

In this section we define and study a cohomology theory specifically
suited for the more general family $x^d + a x$ with $p \nmid d$.
While the growth conditions on the Banach spaces $L(b' ;\rho)$ and
$\c K(b',b; \rho)$ below may seem elaborate, they will allow us to
obtain a detailed description of the associated cohomology, as well
as, provide us with an efficient means of reduction modulo the
operator $D_a$. As an arithmetic application, in \S \ref{S:
FrobEstimates} we are able to explicitly study the action of
Frobenius on the relative cohomology space $\c H_a(b',b)$ defined
below.

Let $d$ be a positive integer relatively prime to $p$. Let $b$ and
$b'$ be two positive real numbers. We will assume throughout this
section that $b \geq b' > 0$. With $\rho \in \bb R$, define the
spaces
\begin{align*}
L(b'; \rho) &:= \{ \sum_{i=0}^\infty B_i a^i : ord(B_i) \geq b'(1
- 1/d)i + \rho, \forall i \geq 0 \} \\
L(b') &:= \bigcup_{\rho \in \bb R} L(b'; \rho) \\
\c K(b', b; \rho) &:= \{ \sum_{i,j \geq 0} B_{ij} a^i x^j :
ord(B_{ij}) \geq  b'(1-1/d)i + bj/d + \rho, \forall i \geq 0 \} \\
\c K(b',b) &:= \bigcup_{\rho \in \bb R} \c K(b',b;\rho) \\
\c V(b', b; \rho) &:= \left( \bigoplus_{i=0}^{d-1} L(b'; \rho)
x^i \right) \cap \c K(b', b; \rho) \\
\c V(b',b) &:= \bigcup_{\rho \in \bb R} \c V(b', b; \rho).
\end{align*}
Notice that $\c K(b',b)$ is an $L(b')$-module.

Fix $\pi \in \bb C_p$ such that $\pi^{p-1} = -p$. Define $D_a$, an
$L(b')$-module endomorphism of $\c K(b',b)$, by
\[
D_a := x \frac{\partial}{\partial x} + \pi(d x^d + ax).
\]
It is useful to keep in mind that, formally,
\begin{equation}\label{E: FormalRelDiff}
D_a = e^{-\pi(x^d+ax)} \circ x \frac{\partial}{\partial x} \circ
e^{\pi(x^d+ax)}.
\end{equation}
Note, we need to use the word ``formally'' because multiplication by
$e^{\pi(x^d+ax)}$ is not an endomorphism of $\c K(b',b)$.

Using this, we may define a \emph{relative cohomology space} as the
$L(b')$-module
\[
\c H_a(b',b) := \c K(b',b) / D_a \c K(b',b).
\]
As the following theorem demonstrates, $\c H_a(b',b)$ is a free
$L(b')$-module of rank $d$.

\begin{theorem}\label{T: RelativeDecomp}
Let $b$ be a positive real number such that $e := b - \frac{1}{p-1}
> 0$. Then, for every $\rho \in \bb R$ we have
\[
\c K(b', b;\rho) = \c V(b', b;\rho) \oplus D_a \c K(b', b;\rho+e)
\]
and
\[
\c K(b', b) = \c V(b', b) \oplus D_a \c K(b', b)
\]
Furthermore, $D_a$ is an injective operator on $\c K(b',b)$.
\end{theorem}

\begin{remark}
We will prove this theorem using Dwork's general method
\cite{DwZetaIHES} which consists of the following six lemmas. Note
that the last three lemmas follow automatically from the validity of
the first three.
\end{remark}

\begin{lemma}\label{L: Decomp_H}
Define $H := \pi(d x^d + ax)$. Then for every $\rho \in \bb R$
\[
\c K(b', b; \rho) = \c V(b', b; \rho) + H \c K(b', b; \rho + e).
\]
\end{lemma}

\begin{proof}
Note, it is sufficient to prove this for $\rho=0$. Let $u =
\sum_{n, s \geq 0} B_{s,n} a^s x^n \in \c K(b', b; 0)$. We may
write this as
\[
u = \sum_{n \geq 0} C_n x^n \quad \text{where} \quad C_n := \sum_{s
\geq 0} B_{s,n} a^s.
\]
For each $j \geq 0$, define
\[
A_j := \sum_{i \geq 0} (-1)^i \left(\frac{a}{d}\right)^i
C_{di+d+(j-i)}.
\]
With $A_{-1} := 0$, a calculation shows that $u = P + H Q$ where
\[
P := \sum_{j=0}^{d-1} ( C_j - \frac{a}{d} A_{j-1}) x^j \quad
\text{and} \quad Q := \frac{1}{d \pi} \sum_{j \geq 0} A_j x^j.
\]
It is not hard to show $P(a,x) \in \c K(b',b;0)$ and $Q(a,x) \in \c
K(b',b;e)$.
\end{proof}

\begin{lemma}
$H \c K(b', b) \cap \c V(b', b) = \{0\}$.
\end{lemma}

\begin{proof}
Suppose
\[
\pi(dx^d + ax) \sum_{j\geq 0} B_j x^j = C_0 + C_1 x + \cdots +
C_{d-1} x^{d-1}
\]
where $\sum_{j \geq 0} B_j x^j \in \c K(b',b)$, $B_j \in
L(b')$, and the right-hand side is in $\c V(b',b)$. For each $j
\geq 0$, the coefficient of $x^{d+j}$ on the left-hand side of
the above equals zero:
\[
\pi d B_j + \pi a B_{d+j-1} = 0 \quad \Longrightarrow\quad B_j =
-\frac{a}{d} B_{d+j-1}.
\]
Iterating this $i$ times yields
\[
B_j = \left( \frac{-a}{d} \right)^i B_{id -i +j}.
\]
Thus, $B_j$ is divisible by $a^i$ for all $i \geq 1$, and so
$B_j$ must equal zero.
\end{proof}

\begin{lemma}
If $f \in \c K(b', b)$ and $H f \in \c K(b', b; \rho)$, then $f \in
\c K(b', b; \rho  + e)$.
\end{lemma}

\begin{proof}
Write $f = \sum_{j \geq 0} A_j x^j$ where $A_j := \sum_{i \geq 0}
G_{ij} a^i$. Now,
\[
H \sum_{j \geq 0} A_j x^j = \sum_{j \geq 0} C_j x^j \in \c K(b', b;
\rho)
\]
where $C_j := \sum_{i \geq 0} B_{ij} a^i$. In particular, this means
\[
\sum_{j \geq d} [d A_{j-d} + a A_{j-1} ] x^j = \frac{1}{\pi} \sum_{j
\geq d} C_j x^j.
\]
Thus,
\[
A_j = -\frac{a}{d} A_{j+d-1} + \frac{1}{d\pi} C_{j+d} \quad \text{
for every } j \geq 0.
\]
Iterating the $A$'s on the right-hand side $N$ times, we obtain the
formula:
\[
A_j = \frac{(-a)^N}{d^N} A_{j + N(d-1)} + \sum_{i=0}^{N-1}
\frac{(-a)^i}{d^i} \frac{1}{d \pi} C_{j + d(i+1)-i}.
\]
Specializing $a$ with $ord(a) + b'(1-\frac{1}{d}) > 0$, then finding
$\delta \in \bb R$ such that $f \in \c K(b',b;\delta)$ we have
\[
ord(\frac{(-a)^N}{d^N} A_{j + N(d-1)}) \geq N[ ord(a) +
b(1-\frac{1}{d})] + bj/d + \delta.
\]
Since $b \geq b'$, the coefficient of $N$ is positive. Thus, for $a$
specialized, letting $N$ tend to infinity we have:
\begin{align*}
A_j &= \sum_{i \geq 0} \left( \frac{-a}{d} \right)^i \frac{1}{d
\pi} C_{j + d(i+1)-i} \\
&= \sum_{r \geq 0} G_{r, j} a^r
\end{align*}
where
\[
G_{r,j} := \sum_{i=0}^r (-1)^i \frac{1}{\pi d^{i+1}} B_{r-i, j +
d(i+1)-i}.
\]
It follows that $f(a,x) = \sum_{j, r \geq 0} G_{r,j} a^r x^j \in \c
K(b',b;\rho + e)$.
\end{proof}

\begin{lemma}\label{L: Decomp_D}
For every $\rho \in \bb R$, $\c K(b', b;\rho) = \c V(b', b; \rho) +
D_a \c K(b', b; \rho + e)$.
\end{lemma}

\begin{proof}
It is sufficient to prove this for $\rho = 0$. Let $f \in \c
K(b',b;0)$. Set $f^{(0)} := f$. Then there exists a unique
$\eta^{(0)} \in \c V(b',b;0)$ and $\xi^{(0)} \in \c K(b',b;e)$ such
that
\[
f^{(0)} = \eta^{(0)} + H \xi^{(0)}.
\]
Define
\begin{align*}
f^{(1)} :&= f^{(0)} - \eta^{(0)} - D_a \xi^{(0)} \\
&= - x \frac{\partial}{\partial x} \xi^{(0)} \in \c K(b',b;e).
\end{align*}
Then, there exists a unique $\eta^{(1)} \in \c V(b',b;e)$ and
$\xi^{(1)} \in \c K(b', b; 2e)$ such that
\[
f^{(1)} = \eta^{(1)} + H \xi^{(1)}.
\]
Define
\[
f^{(2)} := f^{(1)} - \eta^{(1)} - D_a \xi^{(1)}.
\]
Continuing $h$ times we get:
\[
f^{(h)} := f^{(h-1)} - \eta^{(h-1)} - D_a \xi^{(h-1)}.
\]
Adding all these together, we obtain
\[
f^{(h)} = f^{(0)} - \sum_{i=0}^{h-1} \eta^{(i)} - D_a
\sum_{i=0}^{h-1} \xi^{(i)} \in \c K(b',b; he).
\]
Thus, as $h \rightarrow \infty$, $f^{(h)} \rightarrow 0$ in $\c
H(b', b)$, leaving
\[
f = \sum_{i \geq 0} \eta^{(i)} + D_a \sum_{i \geq 0} \xi^{(i)} \in
\c V(b',b) + D_a \c K(b',b).
\]
\end{proof}

\begin{lemma}\label{L: DaLift}
If $f \in \c K(b', b)$ and $D_a f \in \c K(b', b; \rho)$, then $f
\in \c K(b', b; \rho + e)$.
\end{lemma}

\begin{proof}
If $f \ne 0$, then we may choose $c \in \bb R$ such that $f \in \c
K(b',b;c)$ but $f \notin \c K(b', b; c+ e)$. By hypothesis, $D_a f
\in \c K(b',b;\rho)$, thus
\[
H f = D_a f - x \frac{\partial}{\partial x} f \in \c K(b',b;\rho) +
\c K(b', b; c) = \c K(b', b; l)
\]
where $l := \min\{ \rho, c \}$. Thus, $f \in \c K(b',b; l + e)$
which means $l \ne c$. Hence $l = \rho$ as desired.
\end{proof}

\begin{corollary}\label{C: RelativeInjectiveD}
$ker(D_a | \c K(b', b)) = 0$.
\end{corollary}

\begin{proof}
Suppose $D_a(f) = 0$. Since $0 \in \c K(b', b; \rho)$ for all $\rho
\in \bb R$, by the previous lemma, $f \in \c K(b', b; \rho + e)$ for
all $\rho \in \bb R$. However, the only element with this property
is 0.
\end{proof}

\begin{lemma}\label{L: InterDecomp}
$D_a \c K(b', b) \cap \c V(b', b) = \{0\}$.
\end{lemma}

\begin{proof}
Let $f \in D_a \c K(b',b)\cap \c V(b',b)$ and suppose $f \ne 0$.
Choose $\rho \in \bb R$ such that $f \in \c K(b',b; \rho)$ but $f
\notin \c K(b',b;\rho+e)$. Also, let $\eta \in \c K(b',b)$ such that
$D_a \eta = f$. By Lemma \ref{L: DaLift}, $\eta \in \c K(b', b; \rho
+e)$ and so $f - H\eta = x \frac{\partial}{\partial x} \eta \in \c
K(b',b; \rho + e)$. Thus, there exists $\zeta \in \c V(b',b;\rho +
e)$ and $w \in \c K(b',b; \rho + 2e)$ such that
\[
x \frac{\partial}{\partial x} \eta = \zeta + H w.
\]
This means $f = \zeta + H(\eta + w)$. Since $f \in \c V(b',b)$ we
must have $f = \zeta \in \c V(b',b;\rho + e)$. But this contradicts
our choice of $\rho$.
\end{proof}

\begin{remark}\label{R: divisible}
Notice that in the proof of Lemma \ref{L: Decomp_H} that if $u$ is
divisible by $x$, then $P$ is also divisible by $x$ but $Q$ need not
be. From this and the proof of Lemma \ref{L: Decomp_D}, if $f \in \c
K(b',b; \rho)$ is divisible by $x$ then when we write $f =
\sum_{\geq 0} \eta^{(i)} + D_a \sum_{i \geq 0} \xi^{(i)}$ in the
decomposition of $\c K(b',b;\rho) = \c V(b',b;\rho) \oplus D_a \c
K(b',b;\rho+e)$, then $\sum \eta^{(i)}$ is divisible by $x$ but
$\sum \xi^{(i)}$ need not be. This will be important when we define
the relative primitive cohomology in \S \ref{SS: RelPrimCohom}.
\end{remark}

%%% ----------------------------------------------------------------------
\subsection{Relative Dual Cohomology}\label{S: RelDualThy}

In the last section, we saw that $\{ x^i \}_{i=0}^{d-1}$ is a basis
of the $L(b')$-module $\c H_a(b',b) := \c K(b', b)/ D_a \c K(b',b)$.
In this section, we will describe a dual space of $\c H_a(b',b)$ and
compute its dual basis via a nondegenerate pairing.

Let $b^*$ be a positive real number with $b > b^*$. Define the
$L(b')$-module
\[
\c K^*(b',b^*) := \left\{ \sum_{i, j = 0}^\infty B_{ij} a^i x^{-j} |
\inf_{i, j} ( \ord(B_{ij}) - (b'(1 - \frac{1}{d})i - b^*j/d) ) >
-\infty \right\}.
\]
Define the pairing $\langle \cdot, \cdot \rangle:  \c K(b',b) \times
\c K^*(b', b^*) \rightarrow L(b')$ as follows. With $g^* \in \c
K^*(b', b^*)$ and $h \in \c K(b', b)$, thinking of the variable $a$
as a constant, define
\[
\langle h, g^* \rangle := \text{constant term w.r.t. $x$ in the
product } h g^*.
\]
For clarity, we will write this out explicitly: with $g^*(a,x) :=
\sum_{j \geq 0} A_j(a) x^{-j} \in \c H^*(b', b^*)$ and $h(a,x) :=
\sum_{j \geq 0} B_j(a) x^j \in \c K(b', b)$,
\[
\langle h, g^* \rangle := \sum_{s = 0}^\infty G_s a^s \quad
\text{where} \quad G_s := \sum_{r = 0}^\infty \sum_{i + j = s}
A_{ir} B_{jr}.
\]
This defines a \emph{perfect} pairing (or nondegenerate pairing) in
the sense that no nonzero element $g^* \in \c K^*(b', b^*)$ exists
such that $\langle h, g^* \rangle = 0$ for every $h \in \c K(b', b)$
and {\it vice versa}.

Define the truncation operator $Trunc_x: \bb C_p[[a, x^{\pm 1}]]
\rightarrow \bb C_p[[a, x^{-1}]]$ by linearly extending
\[
Trunc_x( x^n) :=
\begin{cases}
x^n & \text{if } n \leq 0 \\
0 & \text{if } n > 0
\end{cases}.
\]
Define the $\c K^*(b', b^*)$-endomorphism
\begin{align}\label{D: D_star}
D_a^* :&= -x \frac{\partial}{\partial x} + Trunc_x[ \pi(d x^d + a x)] \notag \\
&= Trunc_x[ e^{\pi(x^d + a x)} \circ - x\frac{\partial}{\partial x}
\circ e^{-\pi(x^d + a x)}].
\end{align}
The operators $D_a$ and $D_a^*$ are dual to one another with respect
of the pairing. That is, for every $g^* \in \c K^*(b', b^*)$ and $h
\in \c K(b', b)$, we have
\[
\langle h, D_a^* g^* \rangle = \langle D_a h, g^* \rangle.
\]
To see this, notice that the operator $Trunc_x$ does not affect
constant terms (with respect to $x$), and so, $Trunc_x[\pi(d x^d +
ax)]$ is dual to $\pi(d x^d + ax)$. Next, that $-x
\frac{\partial}{\partial x}$ is dual to $x \frac{\partial}{\partial
x}$ follows from the Leibniz identity
\[
x \frac{d}{dx} (h g^*) = x \frac{d}{dx}(h) g^* + h
x\frac{d}{dx}(g^*)
\]
since the left-hand side has no constant term.

Define the $L(b')$-module
\[
\c R_a(b',b^*) := ker( D_a^* | \c K^*(b', b^*)).
\]
We wish to show that $\c R_a(b',b^*)$ is the dual space of $\c
H_a(b',b)$.

\begin{theorem}\label{T: RelDualBasis}
Let $b > \frac{1}{p-1} > b' > 0$, and $b \geq b^* > 0$. Then $\c
R_a(b',b^*)$ is the algebraic dual of $\c H_a(b',b)$. Furthermore,
dual to the basis $\{x^i\}_{i=0}^{d-1}$ of $\c H_a(b',b)$ is the
basis $\{g_i^*(a)\}_{i=0}^{d-1} \subset \c R_a(b',b^*)$ which takes
the explicit form: $g_0^* = 1$, and for $i = 1,2, \ldots, d-1$
\[
g_i^*(a) = \sum_{r=i}^{i+d-2} \sum_{l=0}^\infty \sum_{j=0}^{\left\lfloor
\frac{l-(r-i)}{d} \right\rfloor} \frac{\bb Z}{(d \pi)^{(d-1)j + (r-i)}}
\left( \frac{-a}{d} \right)^{l-(r-i)-dj} x^{-((d-1)l+r)}.
\]
where ``$\bb Z$'' indicates some determinable nonzero integer.
\end{theorem}

\begin{proof}
It is clear that $g_0^* = 1$ is the dual vector of $1 \in \c
H_a(b',b)$. For $i = 1,\ldots, d-1$, let $g_i^*(a) := x^{-i} +
\sum_{j \geq d} B_j^{(i)} x^{-j}$. Let us determine these
$B_j^{(i)}$.

From $D_a^*(g_i^*(a)) = 0$, the coefficients $B_j^{(i)}$ must
satisfy the recurrence relation
\[
B_{d+j}^{(i)} = \frac{-j}{d \pi} B_j^{(i)} - \frac{a}{d}
B_{j+1}^{(i)}.
\]
Thus, from the initial conditions given by requiring $\langle
x^j, g_i^*(a) \rangle = \delta_{ij}$, where $\delta_{ij}$ is
the Kronecker delta symbol, each $B_j^{(i)}$ may be uniquely
solved by a series of the form indicated above.

Next, notice that the $p$-adic order of the coefficient of
$a^{l-(r-i)-dj} x^{-((d-1)l+r)}$ in $g_i^*(a)$ is bounded below by
$-\frac{(d-1)j + (r-i)}{p-1}$ since $d$ is a $p$-adic unit. Since $b
\geq \frac{1}{p-1} \geq b' > 0$, we have
\[
-\frac{(d-1)j + (r-i)}{p-1} \geq b' \left(1-\frac{1}{d} \right)(
l-(r-i)-dj) - \frac{b}{d}((d-1)l+r).
\]
This means $g_i^*(a) \in \c K^*(b', b^*; 0)$.

To show $\{ g_i^*(a) \}_{i=0}^{d-1}$ is a basis of $\c R_a(b',b^*)$,
we need only show that they span $\c R_a(b',b^*)$ since it is clear
that they are linearly independent. This is demonstrated as follows:
Let $h^* \in \c R_a(b',b^*)$ and define $h_i := \langle x^i, h^*
\rangle$ for each $i = 0,1,\ldots, d-1$. Next, define $\tilde h^* :=
h^* - \sum_{i=0}^{d-1} h_i g_i^*(a)$. We wish to show $\tilde h^* =
0$. Now, $\tilde h^* = \sum_{j \geq 0} \tilde h_j^*(a) x^{-j}$, and
by construction, $\tilde h^*_i = 0$ for $i = 0,1,\ldots, d-1$. Using
the reduction formula $d \pi x^{n+d} = D_a(x^n) - n x^n - \pi a
x^{n+1}$ and $D_a^*(\tilde h^*) = 0$, it is easy to see that
\[
\tilde h^*_{n + d} = \langle x^{d+n}, \tilde h^* \rangle =
\sum_{i=0}^{d-1} A_i \tilde h^*_i = 0
\]
where $A_i \in L(b')$. Hence, $\tilde h^* = 0$ as desired.
\end{proof}

%%% ----------------------------------------------------------------------
\subsection{Relative Dwork Operators}\label{S: RelDwOp}

In this section, we define the Dwork operator $\bar \alpha(a)$ on
the relative cohomology $\c H_a(b',b)$, and its dual $\bar
\alpha^*(a)$ on $\c R_a(b',b^*)$. Throughout this section we will
fix real numbers $b$ and $b'$ such that $\frac{p-1}{p} \geq b >
\frac{1}{p-1}$ and $\frac{b}{p} \geq b'
> 0$. Fix $\pi \in \bb C_p$ such that $\pi^{p-1} = -p$.

Dwork's first splitting function on $\bb F_p$ is $\theta(t) :=
e^{\pi(t-t^p)} = \sum_{i=0}^\infty \theta_i t^i$. It is well-known
\cite{pM70} that $ord(\theta_i) \geq \frac{p-1}{p^2} i$ for every
$i$. Also, $\theta_i = \frac{\pi^i}{i!}$ for each $i = 0,1,2,
\ldots, p-1$. Next, define
\[
F(a,x) := \frac{\exp \pi(x^d + ax)}{\exp \pi(x^{dp} + a^p x^p)} =
\theta(x^d) \theta(ax).
\]
Writing $F(a,x) = \sum_{r \geq 0} H_r(a) x^r$ with $H_r(a) :=
\sum_{i=0}^{\lfloor r/d \rfloor} \theta_i \theta_{r-di} a^{r-di}$,
then the coefficient of $a^{r-di} x^r$ satisfies $ord(\theta_i
\theta_{r-di}) \geq  b'(1- \frac{1}{d})(r-di) + \frac{(b/p)}{d} r$.
Consequently, $F(a,x) \in \c K(b',b/p)$. It follows that
multiplication by $F(a,x)$ is an endomorphism of $\c K(b',b/p)$.

From the Cartier operator defined as
\[
\psi_x : \c K(b',b/p) \rightarrow \c K(b',b) \quad \text{takes}
\quad \sum_{i=0}^\infty B_i(a) x^i \longmapsto \sum_{i=0}^\infty
B_{pi}(a) x^i,
\]
we may define the Dwork operator
\[
\alpha(a) := \psi_x \circ F(a,x): \c K(b',b) \rightarrow \c K(b',b).
\]
It is useful to keep in mind the formal identity:
\begin{equation}\label{E: FormalRelFrob}
\alpha(a) = e^{-\pi (x^d + a^p x)} \circ \psi_x \circ e^{\pi (x^d +
ax)}.
\end{equation}

Since $p x \frac{\partial}{\partial x} \circ \psi_x = \psi_x \circ x
\frac{\partial}{\partial x}$, and for any $u \in \bb C_p[[x]]$,
$u(x) \circ \psi_x = \psi_x \circ u(x^p)$, equations (\ref{E:
FormalRelDiff}) and (\ref{E: FormalRelFrob}) demonstrate that
$\alpha(a) \circ D_a = p D_{a^p} \circ \alpha(a)$. Consequently,
$\alpha(a)$ induces two $L(b')$-linear maps on the relative
cohomology: $\bar \alpha(a): \c H_a(b',b) \rightarrow \c
H_{a^p}(b',b)$ and $\bar \alpha(a): ker(D_a | \c K(b',b))
\rightarrow ker(D_{a^p}| \c K(b',b))$.

\begin{remark}\label{R: nondeg}
Let us show $\bar \alpha(a)$ is an isomorphism by creating a
right-inverse on the co-chain level which will reduce to an
isomorphism on cohomology. Denote by $F_{-\pi}(a,x)$ the replacement
of $\pi$ by its conjugate $-\pi$ in $F(a,x)$. Notice that
$F_{-\pi}(a,x) = F(a,x)^{-1} \in \c K(b',b/p)$. Define the Frobenius
map
\[
\Phi_x : \bb C_p[[a, x^{\pm 1}]] \rightarrow \bb C_p[[a, x^{\pm p}]]
\quad \text{by linearly extending} \quad x \mapsto x^p.
\]
Next, define the endomorphism $\alpha(a)^{-1} := F_{-\pi}(a,x) \circ
\Phi_x : \c K(b',b) \rightarrow \c K(b',b)$ and note
\[
\alpha(a)^{-1} = e^{-\pi(x^d+ax)} \circ \Phi_x \circ e^{\pi(x^d+a^p
x)}.
\]
From this and (\ref{E: FormalRelFrob}), $\alpha(a) \circ
\alpha(a)^{-1} = id$. Warning, $\alpha(a)^{-1}$ is \emph{not} a left
inverse of $\alpha(a)$.

Next, from (\ref{E: FormalRelDiff}) and the relation $p \Phi_x \circ
x \frac{\partial}{\partial x} = x \frac{\partial}{\partial x} \circ
\Phi_x$, we see that
\[
p \alpha(a)^{-1} \circ D_{a^p} = D_a \circ \alpha(a)^{-1}.
\]
Hence, $\alpha(a)^{-1}$ induces a mapping $\bar \alpha(a)^{-1}: \c
H_{a^p}(b',b) \rightarrow \c H_a(b',b)$ which satisfies $\bar
\alpha(a) \bar \alpha(a)^{-1} = id$. Since these maps are acting on
free $L(b')$-modules of finite rank, $\bar \alpha(a)$ must be an
isomorphism.
\end{remark}

\medskip\noindent{\bf Dual Dwork Operator.} Define the operator
$\alpha^*(a): \c K^*(b', b^*) \rightarrow \c K^*(b', b^*)$ by the
following:
\begin{align}\label{E: FormalDual}
\alpha^*(a) :&= Trunc_x[ F(a,x) \circ \Phi_x ] \notag \\
&= Trunc_x[ e^{\pi(x^d + ax)} \circ \Phi_x \circ e^{-\pi(x^d + a^p
x)}].
\end{align}

Notice that $\psi_x$ and $\Phi_x$ are dual operators in the sense
that for every $g^* \in \c K^*(b', b^*)$ and $h \in \c K(b', b)$, we
have $\langle \psi_x h, g^* \rangle = \langle h, \Phi_x g^*
\rangle$. Next, since the operator $Trunc_x$ does not kill any
constant terms (with respect to $x$), $Trunc_x[F(a,x)]$ is dual to
$F(a,x)$. Putting these together yields the duality between
$\alpha^*(a)$ and $\alpha(a)$. By taking the dual of $\alpha(a)
\circ D_a = p D_{a^p} \circ \alpha(a)$,  we have $p \alpha^*(a)
\circ D_{a^p}^* = D_a^* \circ \alpha^*(a)$. Therefore, $\alpha^*(a)$
induces a mapping $\alpha^*(a): \c R_{a^p}(b',b^*) \rightarrow \c
R_a(b',b^*)$. Moreover, $\bar \alpha(a): \c H_a(b',b) \rightarrow \c
H_{a^p}(b',b)$ and $\alpha^*(a): \c R_{a^p}(b', b^*) \rightarrow \c
R_a(b',b^*)$ are dual to one another with respect to the pairing
$\langle \cdot, \cdot \rangle$.

%%%-----------------------------------------------------------------------
\subsection{Relative Primitive Cohomology and its Dual}\label{SS:
RelPrimCohom}

Define the $L(b')$-module
\[
x \c K(b',b) := \{ g \in \c K(b',b) | g \text{ is divisible by $x$}
\}.
\]
Notice that $\alpha(a)$ and $D_a$ are well-defined endomorphism of
$x \c K(b',b)$. Define $\c M_a(b',b) := x \c K(b',b) / D_a \c
K(b',b)$. From Theorem \ref{T: RelativeDecomp} and Remark \ref{R:
divisible}, $\c M_a(b',b)$ is a free $L(b')$-module with basis
$\{x^i\}_{i=1}^{d-1}$.

Since we still have $p D_{a^p} \circ \alpha(a) = \alpha(a) \circ
D_a$, $\alpha(a)$ induces a mapping
\[
\bar \alpha(a): \c M_a(b',b) \rightarrow \c M_{a^p}(b',b).
\]
Using a pairing identical to that between $\c H_a(b',b)$ and $\c
R_a(b',b^*)$, we see that the subspace generated by $1$ in $\c
R_a(b',b^*)$ is the annihilator of $\c M_a(b',b)$. Hence, the dual
of $\c M_a(b',b)$ is the free $L(b')$-module $\c R_a'(b',b^*) := \c
R_a(b',b^*) / \langle 1 \rangle$ with basis $\{
g_i^*(a)\}_{i=1}^{d-1}$. Dual to $\bar \alpha(a)$ is $\bar
\alpha^*(a)$, induced by $\alpha^*(a)$, and this map takes $\c
R_{a^p}'(b',b^*)$ bijectively onto $\c R_a'(b',b^*)$. Notice that
$\bar \alpha^*(a)(1) = 1$.

%%%-----------------------------------------------------------------------
\subsection{Relative Functional Equation}\label{S: RelFunEqu}

So far, we have fixed a solution $\pi$ of the equation $x^{p-1} =
-p$. Since we are assuming $p$ is odd, $-\pi$ is also a solution.
Let us see how we would have proceeded had we used $-\pi$. But
before we do this, we need to modify our current notation to keep
things clear. Let us denote the operator $D_a$ by $D_{\pi, a}$, the
space $\c M_a(b',b)$ by $\c M_{\pi, a}(b',b)$, and the Dwork
operator $\bar \alpha(a)$ by $\bar \alpha_\pi(a)$. Also, denote by
$\c R_{\pi,a}(b',b^*)$ and $\alpha_\pi^*(a)$ the dual space $\c
R_a(b',b^*)$ and dual Dwork operator $\alpha^*(a)$.

Now, had we used $D_{-\pi, a} := x\frac{\partial}{\partial x} - \pi
(d x^d + ax)$, then $\c H_{-\pi, a}(b',b) := \c K(b',b) / D_{-\pi,
a} \c K(b',b)$ would still be a free $L(b')$-module with basis
$\{x^i\}_{i=0}^{d-1}$. Let $\c M_{-\pi, a}(b',b) := x \c K(b',b) /
D_{-\pi,a} \c K(b',b)$.

Using the same pairing as that between $\c H_{\pi,a}(b',b)$ and $\c
R_{\pi}(b',b^*)$, the dual of $\c H_{-\pi, a}(b',b)$ is the space
$\c R_{-\pi, a}(b',b^*) := ker( D_{-\pi, a}^* | \c K^*(b',b^*))$.
Clearly, a basis for the latter space is $\{ g_{-\pi, i}^*
\}_{i=0}^{d-1}$, where we have replaced $\pi$ for $-\pi$ in the
formulas given in Theorem \ref{T: RelDualBasis}. Next, dual to
$\alpha_{-\pi}(a) := \psi_x \circ F_{-\pi}(a,x)$ is
$\alpha_{-\pi}^*(a) := Trunc_x \circ F_{-\pi}(a,x) \circ \Phi_x$.
Notice that $\alpha_{-\pi}^*(a): \c R_{-\pi, a^p}(b',b^*)
\rightarrow \c R_{-\pi, a}(b',b^*)$.

Define $\Theta_{-\pi, a} := -x \frac{\partial}{\partial x} -
\pi( d x^d + ax)$; it is useful to view this map as both
$D_{-\pi,a}^*$ without the truncation operator, and as
$-D_{\pi,a}$. By definition, $\c R_{-\pi, a}(b', b^*)$ is the
kernel of $D_{-\pi, a}^* = Trunc_x \circ \Theta_{-\pi, a}$.
Thus, if $\xi^* \in \c R_{-\pi, a}(b',b^*)$ then $\Theta_{-\pi,
a} \xi^*$ will consist of only a finite number of positive
powers of $x$. That is, $\Theta_{-\pi,a}$ defines an
$L(b')$-linear map
\[
\Theta_{-\pi, a}: \c R_{-\pi,a}(b',b^*) \rightarrow xL(b')[x].
\]
We view the right-hand side as a subset of $\c K(b,b')$. Therefore,
by reduction, $\Theta_{-\pi,a}$ induces the mapping
\[
\overline{\Theta}_{-\pi, a}: \c R_{-\pi,a}(b',b^*) \rightarrow \c
H_{\pi, a}(b',b).
\]

Notice that $\overline{\Theta}_{-\pi, a}(1) = - D_{\pi, a}(1) =
0$ in $\c H_{\pi, a}(b',b)$. Also, writing $g_{-\pi, i}^*(a,x)
= x^i + \sum_{j \geq d} B_j^{(i)}(a) x^{-j}$ as in the proof of
Theorem \ref{T: RelDualBasis}, we see that
$\overline{\Theta}_{-\pi, a}(g_{-\pi, i}^*(a,x)) = -d \pi
x^{d-i}$ for each $i = 1, 2, \ldots, d-1$ since $g_{-\pi, i}^*$
is in the kernel of $D_{-\pi,a}^*$ and so all negative powers
of $x$ vanish, including the constant term.

It follows that $\overline{\Theta}_{-\pi, a}$ induces an
$L(b')$-module isomorphism
\[
\overline{\Theta}_{-\pi, a}: \c R_{-\pi, a}'(b',b^*) \rightarrow \c
M_{\pi, a}(b',b).
\]

With $\xi^* \in \c R_{-\pi, a}'(b',b)$, we have
\[
e^{-\pi(x^d +ax)} \circ \Phi_x \circ e^{\pi(x^d+a^p x)}(\xi^*) = \bar
\alpha_{-\pi}^*(a)(\xi^*) + \eta
\]
where $\eta \in L(b')[x]$. Applying $\Theta_{-\pi,a}$ to both sides
we have
\[
p e^{-\pi(x^d +ax)} \circ \Phi_x \circ e^{\pi(x^d+a^p x)} \Theta_{-\pi,
a^p}(\xi^*) = \Theta_{-\pi, a} \bar \alpha_{-\pi}^*(a)(\xi^*) +
\Theta_{-\pi, a}(\eta)
\]
which we may rewrite as
\[
\Theta_{-\pi, a^p}(\xi^*) = p^{-1} \alpha_\pi(a) \circ \Theta_{-\pi,
a} \circ \bar \alpha_{-\pi}^*(a)(\xi^*) - D_{\pi, a^p} \circ
\alpha_\pi(a)(\eta).
\]
It follows that
\begin{equation}\label{E: RelFun}
\overline{\Theta}_{-\pi, a^p} = p^{-1} \bar \alpha_\pi(a) \circ
\overline{\Theta}_{-\pi, a} \circ \bar \alpha_{-\pi}^*(a).
\end{equation}
This is the functional equation. Notice that it relates the
conjugate dual Dwork operator $\bar \alpha_{-\pi}^*(a)$ to the
inverse Dwork operator $\bar \alpha_\pi(a)^{-1}$.

%%% ----------------------------------------------------------------------
\subsection{Frobenius Estimates}\label{S: FrobEstimates}

Let $\f A(a) = (\f A_{ij})$ be the matrix of $\bar \alpha(a)$ with
respect to the basis $\{ x^i \}_{i=1}^{d-1}$. In this section, we
will determine the $p$-adic order of $\f A_{ij}$ as a function of
$a$. Fix $b := (p-1)/p$ and $b' := b/p$. Recall from \S \ref{S:
RelDwOp} that $F(a,x) = \sum_{r \geq 0} H_r(a) x^r \in \c K(b', b';
0)$. Thus, $ord(H_r(a)) \geq b' r /d$ for all $ord(a) > -
\left(\frac{d-1}{d} \right) b'$. For $i = 1,2,\ldots, d-1$, since
$x^i \in \c K(b', b'; -ib'/d)$, we see that $F(a,x) x^i \in \c K(b',
b'; -ib'/d)$, and so $\alpha(a) x^i \in \c K(b', p b'; -ib'/d)$.
From Lemma \ref{L: Decomp_D},
\[
\alpha(a) x^i \subset \sum_{j=1}^{d-1} \f A_{ij}(a) x^j + D_a \c
K(b', b)
\]
for some $\f A_{ij} \in L(b'; \frac{b'}{d}(pj-i))$. This means
\[
ord(\f A_{ij}(a)) \geq \frac{b'(pj-i)}{d} \quad \text{for all} \quad
ord (a) > -\left( \frac{d-1}{d} \right) b'.
\]
A better estimate than this is often needed for applications. For
this, we offer the following exact order for $\f A_{ij}$.

\bigskip\noindent{\bf Notation.} Let $f(a)$ be an analytic function over
$\bb C_p$ convergent on $ord(a) + \rho > 0$. We will write $f(a) =
h(a) + o_\rho(>)$ if there is an analytic function $h(a)$ such that
$ord(f(a)) = ord(h(a))$ for all $ord(a) + \rho > 0$.

\begin{theorem}\label{T: FrobEstimates}
Suppose $d \geq 2$ and $p \geq d+6$. Then there exists $\epsilon >0$
which depends on $p$ and $d$ such that
\[
\f A_{ij}(a) =\frac{\pi^{pj-i-(d-1) r_{ij}}}{r_{ij}!
(pj-i-dr_{ij})!} a^{pj-i-dr_{ij}} + o_{\epsilon}(>)
\]
where $r_{ij} := \left\lfloor \frac{pj-i}{d} \right\rfloor$. (Note:
$\epsilon \rightarrow 0^+$ as $p$ tends to infinity.)
\end{theorem}

\begin{remark}
As a consequence, the $p$-adic absolute values of the entries of the
Frobenius $\f A(a)$ are constant on the unit circle $|a| = 1$.
\end{remark}

\begin{proof}
Recall from \S \ref{S: RelDwOp}, we may write $F(a,x) = \sum_{r \geq
0} H_r(a) x^r$ with $H_r(a) := \sum_{i=0}^{\lfloor r/d \rfloor}
\theta_i \theta_{r-di} a^{r-di} \in L(b'; b' r/d)$. In $x \c
K(b',b)$, we have
\[
\alpha(a) x^i = \psi_x \circ F(a,x) x^i = \sum_{j \geq 1} H_{pj-i}
x^j.
\]
From Lemma \ref{L: Decomp_H}, we have
\[
\alpha(a) x^i = \mu_{i,1} x + \mu_{i,2} x^2 + \cdots \mu_{i,d-1}
x^{d-1} + \pi(x^d + ax) Q,
\]
for some $Q$, where
\begin{equation}\label{E: mu}
\mu_{i,j} = H_{pj-i} - \frac{a}{d} \sum_{r \geq 0} (-1)^r \left(
\frac{a}{d} \right)^r H_{p[dr+d+j-1-r]-i}.
\end{equation}

\medskip\noindent{\bf Step 1:} $H_{pj-i} = \frac{\pi^{pj-i-(d-1) r_{ij}}}{r_{ij}! (pj-i-dr_{ij})!}
a^{pj-i-dr_{ij}} + o_{\epsilon}(>)$ where $r_{ij} := \left\lfloor
\frac{pj-i}{d} \right\rfloor$.

Recall that $\theta_r = \pi^r / r!$ for $r=0,1,\ldots, p-1$.
Consequently, notice that for $\left\lceil \frac{p(j-1)-(i-1)}{d}
\right\rceil \leq r \leq \left\lfloor \frac{pj-i}{d} \right\rfloor$,
since $1 \leq i,j \leq d-1$ and $p > d$, we have both
\[
\theta_r = \frac{\pi^r}{r!} \qquad \text{and} \qquad
\theta_{pj-i-dr} = \frac{\pi^{pj-i-dr}}{(pj-i-dr)!}.
\]
Motivated by this, let us split $H_{pj-i}$ into two sums:
\[
H_{pj-i}(a) = \sum_{r=0}^{\left\lceil \frac{p(j-1)-(i-1)}{d}
\right\rceil - 1} \theta_r \theta_{pj-i-dr} a^{pj-i-dr} + \sum_{r =
\left\lceil \frac{p(j-1)-(i-1)}{d} \right\rceil}^{\left\lfloor
\frac{pj-i}{d} \right\rfloor} \frac{\pi^{pj-i-(d-1)r}}{r!
(pj-i-dr)!} a^{pj-i-dr}.
\]
Step 1 follows after some elementary tedious calculations.

\medskip\noindent{\bf Step 2:} $\mu_{ij} = H_{pj-i} +
o_{\epsilon}(>)$.

Now,
\[
\left( \frac{a}{d} \right)^{r+1} H_{p[dr+d+j-1-r]-i} \in L(b';
(b'/d)[p(dr+d+j-1-r)-i] - b'(r+1)),
\]
and so, the infinite series in (\ref{E: mu}) is an element in $L(b';
(b'/d)(pd-pj-p-i)-b')$. Comparing this with $H_{pj-i}$ from Step 1,
after tedious calculations, Step 2 follows.

\medskip\noindent{\bf Step 3:} $\f A_{ij} = \mu_{ij} +
o_{\epsilon}(>)$.

In the notation of Lemma \ref{L: Decomp_D}, $\mu_{ij} = \eta^{(0)}$,
and so, continuing this notation
\[
\f A_{ij} - \mu_{ij} = \sum_{r \geq 1} \eta^{(r)}
\]
with $\eta^{(r)} \in L(b; (b'/d)(pj-i) + r e)$ where $e := b-
\frac{1}{p-1} > 0$.  Thus, the series on the right-hand side is an
element of $L(b'; (b'/d)(pj-i) + e)$. Step 3 follows.
\end{proof}

%%% ----------------------------------------------------------------------
\section{Fibres}

%%% ----------------------------------------------------------------------
\subsection{$L$-function of the Fibres: $L(x^d + \bar z x, T)$}\label{S:
Fibres}

By elementary Dwork theory \cite{pM70}, Dwork's first splitting
function $\theta(t) := \exp \pi (t - t^p)$ converges on the closed
unit disk $D^+(0,1)$ and $\theta(1)$ is a primitive $p$-th root of
unity in $\bb C_p$. Let $\bar z$ be an element of a fixed algebraic
closure of $\bb F_p$ and let $s := [\bb F_p(\bar z): \bb F_p]$.
Denote by $z$ the Teichm\"uller representative in $\bb C_p$ of $\bar
z$; notice that $z^{p^s-1} = 1$.

Letting $f_{\bar z}(x) := x^d + \bar z x$, we may define for each $n
\in \bb Z_{>0}$ the exponential sum
\[
S_n^*(f_{\bar z}) := \sum_{\bar x \in \bb F_{p^{sn}}^*}
\theta(1)^{Tr_{\bb F_{p^{sn}}/\bb F_p}( \bar x^d + \bar z \bar x)}.
\]
The associated $L$-function is
\[
L^*(f_{\bar z}, T) := \exp( \sum_{n \geq 1} S^*_n(f_{\bar z})
\frac{T^n}{n} ).
\]

Dwork's splitting function $\theta(t)$ defines a $p$-adic analytic
representation of the additive character $\theta(1)^{Tr(\cdot)}$.
With
\[
F(z,x) := \theta(x^d)\theta(zx),
\]
by standard Dwork theory, we have (with $x$ the Teichm\"uller
representative of $\bar x$ in $\bb C_p^*$)
\[
\theta(1)^{Tr_{\bb F_{p^{sn}}/\bb F_p}( \bar x^d + \bar z \bar x)} =
F(z,x) F(z^p, x^p) \cdots F(z^{p^{ns-1}}, x^{p^{ns-1}}).
\]
Therefore,
\[
S_n^*(f_{\bar z}) = \sum_{x \in \bb C_p, x^{p^{ns-1}} = 1} F(z,x)
F(z^p, x^p) \cdots F(z^{p^{ns-1}}, x^{p^{ns-1}}).
\]

We now return to our current situation. Let $\bar z, z$, and $s$ be
as above. Let $b$ and $b'$ be real numbers such that $\frac{p-1}{p}
\geq b > \frac{1}{p-1}$ and $b/p > b' > 0$. Define the spaces
\begin{align*}
K(b)_z &:= (\c K(b',b) \text{ with the variable $a$ specialized at
$z$}) \\
xK(b)_z &:= \{ h \in K(b)_z | h \text{ is divisible by } x \}.
\end{align*}
From Theorem \ref{T: RelativeDecomp} and \S \ref{SS: RelPrimCohom},
\[
H_z := K(b)_z / D_z K(b)_z \quad \text{and} \quad M_z := xK(b)_z /
D_z K(b)_z
\]
are $\bb C_p$-vector space with bases $\{x^i\}_{i=0}^{d-1}$ and
$\{x^i\}_{i=1}^{d-1}$, respectively. Dual to these are the spaces
\begin{align*}
R_z &:= (\c R_a(b',b^*) \text{ with $a$ specialized at $z$})\\
R_z' &:= (\c R_a'(b',b^*) \text{ with $a$ specialized at $z$}).
\end{align*}

Since, with respect to the basis $\{x^i\}_{i=0}^{d-1}$ of both $\c
H_a(b',b)$ and $\c H_{a^p}(b',b)$, the matrix of $\bar \alpha(a)$
has coefficients in $L(b')$, we may specialize $\bar \alpha(a)$ at
$z$. Define $\bar \alpha_z := \bar \alpha(a)|_{a=z}$. Notice, $\bar
\alpha_z: H_z \rightarrow H_{z^p}$ is an isomorphism from Remark
\ref{R: nondeg}.

It is well-known that $\alpha_z$ is a nuclear operator on the space
$K(b)_z$. Thus, since $z^{p^s-1} = 1$, the trace formula for nuclear
operators \cite[Thm 6.11]{pM70} tells us that
\[
(p^s-1)Tr_{nuc}( \alpha_{z,s} | K(b)_z) = p^s Tr(\bar \alpha_{z,s} |
ker(D_z| K(b)_z)) - Tr(\bar \alpha_{z,s} | H_z)
\]
where
\[
\alpha_{z,s} := \alpha_{z^{p^{s-1}}} \circ \cdots \circ \alpha_{z^p}
\circ \alpha_z.
\]
On the other hand, Dwork's trace formula tells us that
\[
(p^s-1)Tr_{nuc}( \alpha_{z,s} | K(b)_z) = \sum_{x^{p^s-1}=1}
F(z^{p^{s-1}}, x^{p^{s-1}}) \cdots F(z^p, x^p) F(z,x).
\]
Observe that the right-hand side is just Dwork's $p$-adic analytic
representation of the character sum $S_1^*(f_{\bar z})$. That is,
\[
(p^s-1)Tr_{nuc}( \alpha_{z,s} | K(b)_z) = S_1^*(f_{\bar z}).
\]
Since $\alpha_{z, sn} = (\alpha_{z,s})^n$ we may use the trace
formula to generalize this to
\[
S_n^*(f_{\bar z}) = p^{sn} Tr(\bar \alpha_{z,s}^n | ker(D_z|
K(b)_z)) - Tr(\bar \alpha_{z,s}^n | H_z).
\]
Using the well-known identity $\exp(-\sum_{n=0}^\infty Tr(A^n)
\frac{T^n}{n}) = det(1 - AT)$ for finite square matrices $A$, we
obtain
\begin{align*}
L^*( f_{\bar z}, T) &= \exp( \sum_{n \geq 1} S^*_n( f_{\bar z}) \frac{T^n}{n} ) \\
&=\frac{det(I - \bar \alpha_{z,s} T | H_z)}{det(I - p^s \bar
\alpha_{z,s} T | ker(D_z | K(b)_z)}.
\end{align*}
From \S \ref{SS: RelPrimCohom}, we know that $\alpha_{z,s}^*$ has 1
as an eigenvector with eigenvalue 1. Thus,
\[
det(1-\alpha_{z,s}^* T | R_z) = (1-T)det(1-\bar \alpha_{z,s}^* T |
R_z').
\]
Also, $ker(D_z|K(b)_z)=0$ by Theorem \ref{T: RelativeDecomp}.
Therefore, since $det(1-\bar \alpha_{z,s}^* T | R_z') = det(1-\bar
\alpha_{z,s} T | M_z)$,
\begin{align*}
L^*(f_{\bar z}, T) &= (1-T) det(1 - \bar \alpha_{z,s} T | M_z) \\
&= (1-T) \prod_{i=1}^{d-1} (1 - \pi_i(\bar z)T) \in \bb
Z[\zeta_p][T].
\end{align*}
Equivalently, for every $n \geq 1$,
\[
-S_n^*(f_{\bar z}) = 1 + \pi_1(\bar z)^n + \cdots + \pi_{d-1}(\bar
z)^n.
\]
Now, had we used $\bar x = 0$ in the definition of the character sum
$S_n^*(f_{\bar z})$, that is,
\[
S_n(f_{\bar z}) := \sum_{\bar x \in \bb F_{p^{sn}}}
\theta(1)^{Tr_{\bb F_{p^{sn}}/\bb F_p}( \bar x^d + \bar z \bar x)},
\]
then since $S_n(f_{\bar z}) = 1 + S_n^*(f_{\bar z})$, we have
\[
L( f_{\bar z}, T) := \exp( \sum_{n \geq 1} S_n( f_{\bar z})
\frac{T^n}{n} ) = det(I - \bar \alpha_{z,s} T | M_z) =
\prod_{i=1}^{d-1} (1 - \pi_i(\bar z)T).
\]
Notice that, from the relative functional equation (\ref{E:
RelFun}), there is a nonzero constant $c$, dependent on $p$, $d$,
and $deg(\bar z)$, such that
\[
c T^{d-1} \overline{L}(f_{\bar z}, p^{-s} T^{-1}) = L(f_{\bar z}, T),
\]
where the bar on the left means complex conjugation. Note, when $d$
is odd, since $-(x^d+ax) = (-x)^d + a(-x)$, the $L$-function has
real coefficients since the exponential sums are in this case real.

%%% ----------------------------------------------------------------------
\bigskip\noindent
{\bf Newton Polygon.} Let us determine the Newton polygon of the
$L$-function $L(f_{\bar z}, T)$. For simplicity, we will assume
$\bar z \in \bb F_p^*$. Set $b := \frac{p-1}{p}$ and $b' := b/p$.
Let $\xi \in \bb C_p$ such that $ord_p(\xi) = b'/d$, and let $\tilde
\xi := \xi^{p-1}$. Consider the basis $\{ \xi^i x^i \}_{i=1}^{d-1}$
for the space $\c M_a(b',b)$. From the beginning of \S \ref{S:
FrobEstimates}, we have
\[
\bar \alpha(a) \xi^i x^i = \sum_{j=1}^{d-1} (\f A_{i,j} \xi^{i-j})
\xi^j x^j \in \c M_{a^p}(b',p b'; 0)
\]
Thus, with respect to this basis, the matrix $\f A(a) = (\f A_{ij})$
of $\bar \alpha(a)$, acting on the right, may be written as
\[
\f A(a) =
\left(%
\begin{array}{cccc}
  A_{0,0} &  \tilde \xi A_{0,1} & \cdots & \xi^{d-1} A_{0, d-1}  \\
  A_{1,0} & \tilde \xi A_{1,0} & \cdots & \xi^{d-1} A_{1,d-1} \\
  \vdots & \vdots & \vdots & \vdots \\
  A_{d-1, 0} & \tilde \xi A_{d-1,1} & \cdots  & \tilde \xi^{d-1} A_{d-1, d-1}\\
\end{array}%
\right)
\]
where $\tilde \xi^j A_{i,j} = \f A_{ij} \xi^{i-j}$ and $ord(\f
A_{i,j}) \geq 0$. Now, $det(1 - \f A(a) T) = \sum_{m=0}^{d-1} c_m
T^m$ where
\[
c_m := (-1)^m \sum_{1 \leq u_1 < u_2 < \cdots < u_m \leq d-1}
\sum_{\sigma} sgn(\sigma) \left( \tilde \xi^{u_1} A_{\sigma(u_1),
u_1} \right) \left( \tilde \xi^{u_2} A_{\sigma(u_2), u_2} \right)
\cdots \left( \tilde \xi^{u_m} A_{\sigma(u_m), u_m} \right)
\]
where the second sum runs over all permutations  $\sigma$ of the
$u_1, u_2, \ldots, u_m$ and $sgn(\sigma)$ is the signature of the
permutation. Our goal is to determine the precise $p$-adic order of
these $c_m$ when $ord_p(a) = 0$. We do this as follows. Notice that
\begin{align*}
c_m &\equiv (-1)^m \sum_{\sigma \in S_m} sgn(\sigma) \left( \tilde
\xi A_{\sigma(1), 1} \right) \left( \tilde \xi^2 A_{\sigma(2), 2}
\right) \cdots \left( \tilde \xi^{m} A_{\sigma(m),
m} \right) \quad mod(\tilde \xi^{m+1}) \\
&= (-1)^m \sum_{\sigma \in S_m} sgn(\sigma) \f A_{\sigma(1),1} \f
A_{\sigma(2),2} \cdots \f A_{\sigma(m),m}
\end{align*}
where $S_m$ is the symmetric group on $\{1, 2, \ldots, m\}$. From
Theorem \ref{T: FrobEstimates}, we have
\[
\f A_{\sigma(1),1} \f A_{\sigma(2),2} \cdots \f A_{\sigma(m),m} =
\frac{\pi^{(p-1)m(m+1)/2 - (d-1) \sum_{j=1}^m r_{\sigma(j),j}}}{
\left( \prod_{j=1}^m r_{\sigma(j),j}! \right) \left( \prod_{j=1}^m
[pj-\sigma(j)- d r_{\sigma(j),j}]! \right) }
a^{(p-1)\frac{m(m+1)}{2} - d \sum_{j=1}^m r_{\sigma(j),j}} +
o_\epsilon(>)
\]

Let $\tau_p$ be the permutation on $\{1, 2, \ldots, d-1\}$ defined
by multiplication by $p$ modulo $d$. Define $\tilde r_j := (pj -
\tau_p(j))/d$. Since
\[
\sum_{j=1}^m r_{\sigma(j), j} = \sum_{j=1}^m \left\lfloor
\frac{\tau_p(j) - \sigma(j)}{d} + \tilde r_j \right\rfloor,
\]
as $\sigma$ runs over the permutations in $S_m$, all of these sums
with $\sigma(j) \leq \tau_p(j)$ for $j =1, \ldots, m$ will be equal
and the largest; all other $\sigma$ will give a strictly smaller
sum. Thus,
\[
c_m = u_m \pi^{(p-1)m(m+1)/2 - (d-1)(\tilde r_1 + \tilde r_2 +
\cdots + \tilde r_m)} a^{(p-1)m(m+1)/2 - d (\tilde r_1 + \cdots +
\tilde r_m)} + o_\epsilon(>)
\]
where
\begin{align*}
u_m :&= (-1)^m \frac{1! 2! \cdots m!}{\tilde r_1! \cdots \tilde r_m!
\tau_p(1)! \tau_p(2)! \cdots \tau_p(m)!} \sum_{\sigma \in S_m,
\sigma \leq \tau_p} sgn(\sigma) \binom{\tau_p(1)}{\sigma(1)} \cdots
\binom{\tau_p(m)}{\sigma(m)} \\
&= (-1)^{m+1} \frac{1}{\tilde r_1! \cdots \tilde r_m!} \prod_{1 \leq r < s \leq m} (\tau_p(s) - \tau_p(r)),
\end{align*}
where the second equality comes from the following identity:

\begin{lemma}[J. Zhu \cite{ZhuAiry}]\label{L: BinSUM}
Fix a positive integer $d$. For every $\tau \in S_d$ and $1
\leq m \leq d$,
\begin{equation}\label{E: binomSUM}
\sum_{\sigma \in S_m, \sigma \leq \tau} sgn(\sigma)
\binom{\tau(1)}{\sigma(1)} \cdots \binom{\tau(m)}{\sigma(m)} = - \frac{1}{1! 2! \cdots m!}
\left(\prod_{j=1}^m \tau(j) \right)
\left(\prod_{1 \leq r < s \leq m} (\tau(s) - \tau(r))\right).
\end{equation}
\end{lemma}

\begin{proof}
First, rewrite (\ref{E: binomSUM}) as
\begin{align*}
\frac{1}{1! 2! \cdots m!} &\sum_{\sigma \in S_m, \sigma \leq \tau}
sgn(\sigma) \prod_{i=1}^m \tau(i) (\tau(i)-1) \cdots (\tau(i)-\sigma(i)+1) \\
&=\frac{1}{1! 2! \cdots m!} \sum_{\sigma \in S_m}
sgn(\sigma) \prod_{i=1}^m \tau(i) (\tau(i)-1) \cdots (\tau(i)-\sigma(i)+1) \\
&= \frac{1}{1! 2! \cdots m!} \det M
\end{align*}
where $M$ is the $m \times m$ matrix
\[
M := \left(
  \begin{array}{ccc}
    \tau(1) & \tau(1)(\tau(1)-1) & \cdots \\
    \tau(2) & \tau(2)(\tau(2)-1) & \cdots \\
              &       \vdots           &         \\
    \tau(m) & \tau(m)(\tau(m)-1) & \cdots \\
  \end{array}
\right).
\]
Now,
\[
det(M) = (-1)^{m+1} det
\left(
  \begin{array}{cccc}
    1 & \tau(1) & \tau(1)(\tau(1)-1) & \cdots \\
    1 & \tau(2) & \tau(2)(\tau(2)-1) & \cdots \\
    \vdots &          &       \vdots           &         \\
    1& \tau(m) & \tau(m)(\tau(m)-1) & \cdots \\
    1 & 0 & 0 & \cdots
  \end{array}
\right).
\]
This matrix may be transformed into a Vandermonde matrix by
column operations:
\begin{align*}
det(M) &= (-1)^{m+1}
det\left(
  \begin{array}{cccc}
    1 & \tau(1) & \tau(1)^2 & \cdots \\
    1 & \tau(2) & \tau(2)^2 &  \\
      &           &  \vdots & \cdots \\
    1 & \tau(m) & \tau(m)^2 & \cdots \\
    1 & 0 & 0 & \cdots
  \end{array}
\right) \\
&= - \left(\prod_{j=1}^m \tau(j) \right) \left(\prod_{1
\leq r < s \leq m} (\tau(s) - \tau(r))\right)
\end{align*}
\end{proof}

Since $u_m$ is a $p$-adic unit, we have
\[
ord_p(c_m) = \frac{m(m+1)}{2} - \left( \frac{d-1}{p-1} \right)
(\tilde r_1 + \tilde r_2 + \cdots + \tilde r_m ).
\]
This proves the following theorem. (Note, the following theorem
was first proven by Zhu \cite{ZhuAiry} for the family $x^d +
ax$ where $a \in \bb Q$ and assuming $p$ is sufficiently
large.)

\begin{theorem}\label{T: LfunctionNP}
If $p \geq d+6$, and $\bar z \in \bb F_p^*$ then writing
\[
L(x^d+ \bar z x, T) = (1 - \pi_1(\bar z) T) \cdots (1-
\pi_{d-1}(\bar z)T)
\]
we have
\[
ord_p(\pi_j(\bar z)) = j - \left( \frac{d-1}{p-1} \right)
\left( \frac{pj - \tau_p(j)}{d} \right)
\]
where $\tau_p$ is the permutation on $\{1,2,\ldots, d-1\}$ defined
by multiplication by $p$ modulo $d$.
\end{theorem}

\begin{figure}
\begin{center}
\includegraphics[scale=0.75]{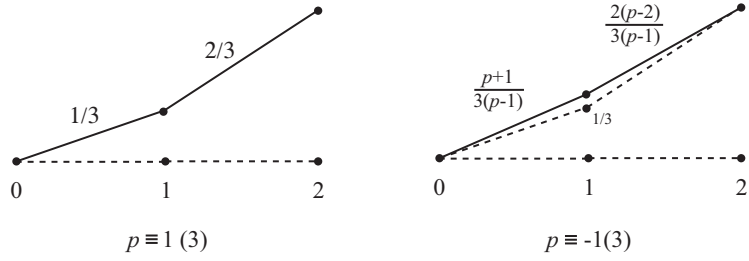}
\caption{Newton polygon of $L(x^3 + \lambda x / \bb F_p,
T)$.}\label{F: AiryNP}
\end{center}
\end{figure}

When $d=3$, Theorem \ref{T: LfunctionNP} says if $p \equiv 1$
mod $3$ then the reciprocal roots $\pi_1(\lambda)$ and
$\pi_2(\lambda)$ may be ordered such that $ord_p \pi_1(\lambda)
= 1/3$ and $ord_p \pi_2(\lambda) = 2/3$ for all $\lambda$. If
$p \equiv -1$ mod $3$, then $ord_p \pi_1(\lambda) =
\frac{p+1}{3(p-1)}$ and $ord_p \pi_2(\lambda) =
\frac{2(p-2)}{3(p-1)}$. See Figure \ref{F: AiryNP}. This was
first proven by Sperber \cite{SperberCubic}.

When $p \equiv -1$ mod($d$) and $p$ is greater than
approximately $2^{(d-1)/2} \cdot d$, Yang \cite{Yang} has
proven Theorem \ref{T: LfunctionNP} by a different method. His
proof is interesting in that he computes the Frobenius over the
chain complex rather than passing to cohomology. The advantage
of this is that the entries of the Frobenius matrix are given
explicitly; the disadvantage is that the Frobenius matrix has
infinitely many rows and columns. Yang's result follows from a
careful diagonalization procedure of this infinite matrix.

\begin{remark}
Blache and F\'erard \cite{BlacheFerard} have given a lower bound for
the generic Newton polygon of the $L$-function of polynomials of
degree $d$ with $p > 3d$. More precisely, in the space of
polynomials of degree $d$, there is a Zariski open set (the
complement of an associated Hasse polynomial) of which any
polynomial lying in this open set will have the Newton polygon of
the associated $L$-function coinciding with this lower bound. For
one-variable polynomials, see also \cite{ZhuAsympVariation},
\cite{ZhuBlacheFerard}. For higher dimensions, see
\cite{WanVariationNP}.
\end{remark}

\begin{remark}
For $d=3$, when $p \equiv 1$ mod $3$ one may show the existence
of a Tate-Deligne mapping (excellent lifting) of the root
$\pi_1(\lambda)$. We suspect no such lifting exists when $p
\equiv -1$ mod $3$.
\end{remark}

%%% ----------------------------------------------------------------------
\subsection{(Appendix) Elementary Entire}\label{S: FrobDecomp}

This section is independent of the rest of the paper. Throughout
this section, we will assume $\bar z^p = \bar z$; the general case
is handled similarly. Notice that
\[
det(1 - \alpha_z T | K(b)_z) = L^*(f_{\bar z}, T) det(1-p \alpha_z
T| K(b)_z).
\]
Using this equation recursively gives
\[
det(1-\alpha T| K(b)_z) = \prod_{i=0}^\infty L^*(f_{\bar z}, p^i T).
\]
That is, $det(1-\alpha T| K(b)_z)$ is the infinite product of the
same polynomial, each scaled by a factor of increasing $p$. We call
any entire function with such a factorization \emph{elementary}. We
wish to show that this infinite product comes from a vector space
decomposition of $K(b)_z$ into an infinite number of finite
dimensional subspaces, each of the same dimension and each related
to one another by the operator $D_z$.

From Theorem \ref{T: RelativeDecomp}, we may write
\begin{equation}\label{E: FibreDecomp}
K(b)_z = V_z \oplus D_z K(b)_z
\end{equation}
where $V_z := \c V(b',b)|_{a=\bar z}$. Notice that $D_z$ is an
endomorphism of $K(b)_z / D_z K(b)_z$, and so, denote its image in
the quotient by $D_z V_z$. Since $ker(D_z|K(b)_z) = 0$, equation
(\ref{E: FibreDecomp}) becomes a recursive equation, which means we
may write
\[
K(b)_z = V_z \oplus D_z V_z \oplus D_z^2 K(b)_z.
\]
Next, since
\[
\begin{CD}
ker(D_z^2|K(b)_z) @>\text{inj}>> K(b)_z @>D_z^2>> K(b)_z
@>\text{surj}>>
\frac{K(b)_z}{D_z^2 K(b)_z}\\
@V p^2 \bar \alpha_z VV @V p^2 \alpha_z VV @VV\alpha_z V @VV\bar
\alpha_z V \\
ker(D_z^2|K(b)_z) @>\text{inj}>> K(b)_z @>D_z^2>> K(b)_z
@>\text{surj}>> \frac{K(b)_z}{D_z^2 K(b)_z}
\end{CD}
\]
where ``inj'' and ``surj'' indicate that the maps are either
injective or surjective, we get
\[
det(1-\alpha_z T| K(b)_z) = det(1- \bar \alpha_z T | K(b)_z / D_z^2
K(b)_z) det(1-p^2 \alpha_z T| K(b)_z).
\]
Since $K(b)_z/D_z^2 K(b)_z$ is isomorphic to $V_z \oplus D_z V_z$,
\[
det(1- \bar \alpha_z T | K(b)_z/D_z^2 K(b)_z) = det(1- \bar \alpha_z
T|V_z) det(1-\bar \alpha_z T| D_z V_z).
\]
It follows that
\[
det(1-\bar \alpha_z T| D_z V_z) = det(1 - p \bar \alpha_z T | V_z) =
L^*(f_{\bar z}, pT).
\]
Of course, this generalizes to
\[
det(1-\bar \alpha_z T| D_z^s V_z) = det(1 - p^s \bar \alpha_z T |
V_z) = L^*(f_{\bar z}, p^s T).
\]
That is, the Fredholm determinant of the operator $\alpha_z$ on the
decomposition
\[
K(b)_z = V_z \oplus D_z V_z \oplus D_z^2 V_z \oplus \cdots
\]
is
\[
det(1- \alpha_z T| K(b)_z) = \prod_{i=0}^\infty det(1-\bar \alpha_z
T | D_z^i V_z) = \prod_{i=0}^\infty L^*(f_{\bar z}, p^i T).
\]
Notice that we have also proven
\[
det(1-\alpha_z T| K(b)_z) = \lim_{s\rightarrow \infty} det(1- \bar
\alpha_z T | K(b)_z / D_z^s K(b)_z).
\]
A dual statement of this takes the form
\[
det(1 - \alpha_z^* T | K(b)_z^*) = \lim_{s \rightarrow \infty} det(1
- \bar \alpha_z^* T | ker((D^*_z)^s | K(b)_z^*)).
\]

%%% ----------------------------------------------------------------------
\section{Variation of Cohomology (Deformation Theory)}\label{S:
DeformTheory}

From now on we will fix $d = 3$. With $\bar z \in \overline{\bb
F}_p$ and $z$ its Teichm\"uller representative in $\bb C_p^*$ define
the spaces
\begin{align*}
K(b)_z^* &:= (\c K^*(b',b) \text{ with $a$ specialized at $z$})
\quad \\
R_z &:= ker(D_z^*| K(b)_z^*) \\
R_z' &:= R_z / \langle 1 \rangle.
\end{align*}
Notice that $R_z$ and $R_z'$ are just $\c R_a(b',b^*)$ and $\c
R_a'(b',b^*)$ with $a$ specialized at $z$, respectively. Through use
of the pairing defined in \S \ref{S: RelDualThy} these three spaces
are algebraically dual to $K(b)_z$, $H_z$, and $M_z$, respectively.

We wish to study how the space $R_z'$ varies as $z$ moves around in
$\bb C_p$. To do this we will define an isomorphism $T_{z,a}$ from
$R_z'$ to $R_a'$ as follows: first, for each $a, z \in \bb C_p$ with
$|a-z| < p^{-b'/3}$ define the isomorphism $T_{z,a}: K(b)_z^*
\rightarrow K(b)_a^*$ by
\[
T_{z,a} := Trunc_x \frac{\exp \pi(x^3+ax)}{\exp \pi (x^3 + zx)} =
Trunc_x \circ e^{\pi(a-z)x}.
\]
Next, using (\ref{D: D_star}),
\[
D_a^* \circ T_{z,a} = T_{z,a} \circ D_z^*.
\]
Also, notice that $T_{z,a}(1) = 1$. Consequently, for any $a \in \bb
C_p$ close enough to $z$ the mapping $T_{z,a}$ induces an
isomorphism $T_{z,a}: R_z' \rightarrow R_a'$. It is important to
notice that $T_{z,z} = I$.

%%% ----------------------------------------------------------------------
\bigskip\noindent{\bf Deformation Theory.}
Using notation from Theorem \ref{T: RelDualBasis}, fix the bases $\{
g_1^*(z), g_2^*(z) \}$ and $\{ g_1^*(a), g_2^*(a) \}$ of $R_z'$ and
$R_a'$, respectively. Observe that each $g_j^*(z)$ is simply a power
series in the variable $x$ with coefficients in $\bb C_p$. Let
$C(z,a) = (c_{i,j})$ be the matrix representation of $T_{z,a}$ with
respect to these bases. Note, we are thinking of $R_z'$ as a column
space, and so, $C(z,a)$ acts on the left.

The pairing $\langle \cdot, \cdot \rangle: M_z \times R_z'
\rightarrow \bb C_p$ allows us to focus on individual entries of
$C(z,a)$:
\begin{equation}\label{E: Cij}
c_{ij} = \langle x^i, T_{z,a}(g_j^*(z)) \rangle = \text{constant
term of $x$ in } x^i e^{\pi(a-z)x} g_j^*(z).
\end{equation}

We wish to demonstrate that $C(z,a)$ satisfies a differential
equation. Differentiating (\ref{E: Cij}) with respect to the
variable $a$, we obtain
\begin{equation}
\frac{d c_{i,j}}{da} = \text{constant term of $x$ in } \pi x^{i+1}
e^{\pi(a-z)x} g_j^*(z) = \langle \pi x^{i+1}, T_{z,a}(g_j^*(z))
\rangle. \notag
\end{equation}
The pairing in this last line is between $M_a$ and $R_a'$. Since
$M_a$ is a $\bb C_p$-vector space spanned by the vectors $x, x^2$ we
see that $\pi x^{i+1}$ is just a scalar multiple of a basis vector
if $i=1$. If $i = 2$ then we need to rewrite $\pi x^3$ in terms of
the basis $\{x, x^2\}$. This is done as follows: using $D_a := x
\frac{\partial}{\partial x} + \pi(3 x^3 + ax)$, we have $D_a(1) =
3\pi x^3 + \pi ax$. Hence, $\pi x^3 \equiv \frac{-\pi a}{3}x$ in
$M_a$.

Putting this together, we have
\begin{equation}\label{E: PF}
\frac{d}{da} C(z,a) = B(a) C(z,a) \quad \text{and} \quad C(z,z) = I
\end{equation}
where $B$ is the $2\times 2$ matrix
\[
B(a) := \left(
  \begin{array}{cc}
    0 & \pi \\
    \frac{-\pi a}{3} & 0 \\
  \end{array}
\right)
\]
In particular, if $(C_1, C_2)^t$ is a solution of this differential
equation, then $C_1$ satisfies the Airy equation $y'' + \frac{\pi^2
a}{3} y = 0$ and $C_2 = \pi^{-1} C_1'$.

%%% ----------------------------------------------------------------------
\bigskip\noindent{\bf Connection with Frobenius.}
From (\ref{E: FormalDual}),
\begin{equation}\label{E: deformFrob}
T_{z,a} \circ \bar \alpha_z^* = \bar \alpha_a^* \circ T_{z^p, a^p}.
\end{equation}
Let us use this and the deformation equation (\ref{E: PF}) to
compute $det(\bar \alpha_a^*)$.

Fix $|z| < p^{2b'/3}$ and let $a \in \bb C_p$ be such that $|a-z| <
p^{-b'/3}$. Define $w_z(a) := det(T_{z,a})$. Notice that $w_z(a)$ is
nonzero on $D^-(z, p^{-b'/3})$ since $T_{z,a}$ is invertible. From
(\ref{E: deformFrob}) we may write
\begin{equation}\label{E: detFrob}
det(\bar \alpha_a^*) =  w_z(a) det(\bar \alpha_z^*)
w_{z^p}(a^p)^{-1}.
\end{equation}
Now, a standard fact from differential equations states that
\[
\frac{d w_z(a)}{da} = Tr(B(a)) w_z(a).
\]
Since $Tr(B(a)) = 0$, $w_z(a)$ is a locally constant function
centered at $z$.

In particular, if we let $z=0$ then $w_0(a)$ is a constant function
on $|a| < p^{-b'/3}$. Since $w_0(0) = 1$, we must have $w_0(a) = 1$
for all $|a| < p^{-b'/3}$. This means $det(\bar \alpha_a^*) =
det(\bar \alpha_0^*)$ for all $|a| < p^{-b'/3}$. However, from \S
\ref{S: RelDwOp}, $\bar \alpha(a)^*$ is an $L(b')$-module
isomorphism, and so, $det(\bar \alpha(a)) \in L(b')^*$. Thus, the
domain of the equality of (\ref{E: detFrob}) may be extended to that
of $L(b')$:
\[
det(\bar \alpha_a^*) = det(\bar \alpha_0^*) \qquad \text{for all }
|a| < p^{2b'/3}.
\]
In other words, $det(\bar \alpha(a)^*)$ is constant for all $a$ in
this domain.

Using either the classical theory of Gauss sums or Dwork theory
\cite{bD69}, it is not hard to calculate the value of $det(\bar
\alpha_0^*)$. However, to describe it, we must first recall the
definition of Gauss sum. Define Dwork's first splitting function
$\theta_s(t) := \exp \pi (t - t^{p^s})$ on the field $\bb F_{p^s}$.
Next, for each $j \in \bb Z$ and $s \in \bb Z_{>0}$ define the Gauss
sum
\[
g_s(j) := - \sum_{t \in \bb C_p, t^{p^s-1} = 1} t^{-j} \theta_s(t)
\]
Let $q := p^s$. Then
\[
det(\bar \alpha_{a,s}^*) =
\begin{cases}
q & \text{if } q \equiv 1 \text{ mod}(3) \\
-g_2(\frac{q^2-1}{3}) & \text{if } q \equiv -1 \text{ mod}(3)
\end{cases}.
\]

%%% ----------------------------------------------------------------------
\bigskip\noindent{\bf Frobenius Action on Solutions of the Differential Equation.}
Let $\f M_z$ denote the field of meromorphic functions (in the
variable $a$) near $z$. Define $\partial_a^* := Trunc_a( - a
\frac{d}{da} + \pi a x)$. Notice that we have the commutative
diagram
\[
\begin{CD}
R_z' \otimes \f M_z @>T_{z,a}>> R_a' \otimes \f M_z \\
@V-a \frac{d}{da}VV  @VV\partial_a^*V \\
R_z' \otimes \f M_z @>T_{z,a}>> R_a' \otimes \f M_z.
\end{CD}
\]
As a consequence we have the following important observation:
$g^*(a) := C_1(a) g_1^*(a) + C_2(a) g_2^*(a) \in R_a' \otimes \f
M_z$ is the image of an element in $R_z' \otimes \f M_z$ by
$T_{z,a}$ which is independent of the variable $a$ if and only if
$\partial_a^*(g^*) = 0$. Notice that this last equation means the
vector $\binom{C_1}{C_2}$ satisfies (\ref{E: PF}).

We wish to define a Frobenius-type endomorphism on the vector space
of local solution at $z$ as follows. Let $\binom{C_1}{C_2}$ be a
local analytic solution of (\ref{E: PF}) about $z$. Thus, $C_1$ and
$C_2$ are power series in $(a-z)$. From the previous paragraph,
there exists $\xi_z^* \in R_z' \otimes \f M_z$, independent of $a$,
such that $T_{z,a}(\xi_z^*) = C_1(a) g_1^*(a) + C_2(a) g_2^*(a)$.
Now, for some constants $b_1$ and $b_2$, $\xi_z^* = b_1 g_1^*(z) +
b_2 g_2^*(z)$. It follows that, if we let $\phi$ denote the
operation of replacing $(a-z)$ by $(a^p - z^p)$, then we have
\[
T_{z^p, a^p}(\xi_{z^p}^*) = C_1^\phi g_1^*(a^p) + C_2^\phi
g_2^*(a^p)
\]
with $\xi_{z^p}^* = b_1 g_1^*(z^p) + b_2 g_2^*(z^p)$. Next, apply
$\bar \alpha_a^*$ to both sides:
\[
\bar \alpha_a^*(C_1^\phi g_1^*(a^p) + C_2^\phi g_2^*(a^p)) = \bar
\alpha_a^* \circ T_{z^p, a^p}(\xi_{z^p}^*) = T_{z,a} \circ \bar
\alpha_z^*(\xi_{z^p}^*).
\]
This equation demonstrates that the left-hand side comes from an
element of $R_z' \otimes \f M_z$ independent of $a$, and so, it
satisfies the differential equation. In particular, if $\f Y$ is a
(local) fundamental solution matrix of (\ref{E: PF}) near $z$, then
there is a constant matrix $M$ such that
\[
\f A_a \f Y^\phi = \f Y M
\]
where $\f A_a$ is the matrix of $\bar \alpha_a^*$. Note, the
matrix $M$ depends on $z$.

The determination of $M$ is an interesting topic. For $z = 0$,
$M = \f A_0$ which is easily calculated via Gauss sums. In the
next section, we will determine $M$ when $1 < |z| < p^{2b'/3}$.

%%% ----------------------------------------------------------------------
\subsection{Behavior near Infinity}\label{S: AiryAsymp}

In this section, we wish to determine the matrix $M$ coming from the
Frobenius action on the set of local solutions of (\ref{E: PF}) at
$z$ near infinity. Once this is done, we will study the solutions of
the $k$-th symmetric power of the scalar equation of (\ref{E: PF}).

With $\pi \in \bb C_p$ such that $\pi^{p-1} = -p$, define the
$p$-adic Airy equation
\begin{equation}\label{E: Airy}
y'' + \frac{\pi^2 a}{3} y = 0.
\end{equation}
Observe that this equation is regular everywhere except at infinity
where it has an irregular singular point. Near infinity, the
asymptotic expansions of (\ref{E: Airy}) will be power series in the
ramified variable $\sqrt{a}$, so, we are naturally led to a change
of variables $a \mapsto a^2$. This changes (\ref{E: Airy}) into
\begin{equation}\label{E: Airy2}
a \tilde y'' - \tilde y' + \frac{4\pi a^5}{3} \tilde y = 0.
\end{equation}
For convenience, let us move infinity to zero by the change of
variable $a \mapsto 1/a$ so that (\ref{E: Airy2}) becomes
\[
a^8 y'' + 3 a^7 y' + \frac{4 \pi^2}{3} y = 0.
\]
To remove the irregular singular point at the origin, we consider
solutions of this differential equation of the form
\[
y(a) = a^{1/2} e^{\kappa \pi a^{-3}} v(a) \quad \text{where} \quad
\kappa := \frac{2 i}{3 \sqrt{3}}.
\]
and $i$ is a fixed square root of $-1$. This means $v(a)$ must
satisfy the differential equation
\[
v'' + (4 a^{-1} - 6 \kappa \pi a^{-4}) v' + (5/4) a^{-2} v = 0.
\]
Since this equation has only a regular singular point at the origin,
we may explicitly solve it using the method of Frobenius. It follows
that a local solution of (\ref{E: Airy2}) about infinity is of the
form
\[
y_1(a) := a^{-1/2} \exp(\kappa \pi a^3) v(1/a)
\]
where
\[
v(a) := \sum_{n=0}^\infty \frac{ \left( \frac{7}{6} \right)_n \left(
\frac{17}{12} \right)_n}{2^n \kappa^n \pi^n (n+1)!} a^{3n}.
\]
Note, $(c)_n := c (c+1)(c+2) \cdots (c+n-1)$. Replacing $i$ with
$-i$ obtains another local solution linearly independent over $\bb
C_p$ from $y_1(a)$:
\[
y_2(a) := a^{-1/2} \exp(-\kappa \pi a^3) \bar v(1/a)
\]
where $\bar v(a)$ is the series defined by replacing $i$ with $-i$
in the coefficients of $v(a)$. Using \cite[p.243]{bD80}, it is not
hard to show that $v(a)$ and $\bar v(a)$ converge on the open unit
disk $D^-(0,1)$.

\begin{remark}\label{R: branch}
Notice that we can center the solutions about any $|z| > 1$ by
choosing a branch of $a^{-1/2}$ and by using $\exp( \kappa \pi
(a-z)^3$. In this case, $v$ and $\bar v$ do not change.
\end{remark}

With $a \mapsto a^2$, the deformation equation (\ref{E: PF}) becomes
\[
\frac{d}{da}C(a) =
\left(%
\begin{array}{cc}
  0 & 2\pi a \\
  -\frac{2\pi a^3}{3} & 0 \\
\end{array}%
\right) C(a).
\]
Thus, if the vector $\binom{C_1}{C_2}$ satisfies this, then $C_1$
satisfies (\ref{E: Airy2}).

For $|z| > 1$, the (local) fundamental solution matrix takes the
form
\[
\f Y(a) = a^{-1/2} V(a) S(a)
\]
where $v_1(a) := v(1/a)$ and $v_2(a) := \bar v(1/a)$, and
\[
V(a) :=
\left(%
\begin{array}{cc}
  v_1(a) & v_2(a) \\
  (-\frac{1}{2}a^{-2} + 3 a^2 \kappa \pi)v_1 + a^{-1} v_1'(a)  & (-\frac{1}{2}a^{-2} - 3 a^2 \kappa \pi)v_2 + a^{-1} v_2'(a) \\
\end{array}%
\right)
\]
and
\[
S(a) :=
\left(%
\begin{array}{cc}
  e^{\kappa \pi a^3} & 0 \\
  0 & e^{-\kappa \pi a^3} \\
\end{array}%
\right).
\]
From the previous section, there is a constant matrix $M := \left(%
\begin{array}{cc}
  m_1 & m_2 \\
  m_3 & m_4 \\
\end{array}%
\right)$ such that $\f A_{a^2} \f Y^\phi = \f Y M$. That is,
\begin{equation}\label{E: DefFrob}
\f A_{a^2} a^{-p/2} V(a^p) S(a^p) = a^{-1/2} V(a) S(a) M.
\end{equation}
Let us view this equation over the ring $\hatR$ of analytic
functions on $1 < |a| < e$, $e > 1$ unspecified, adjoined the
algebraic (over $\hatR$) functions $\sqrt{a}$, $e^{\kappa \pi a^3}$,
$e^{-\kappa \pi a^3}$, $e^{\kappa \pi a^{3p}}$, $e^{-\kappa \pi
a^{3p}}$ satisfying the compatibility relations: $e^{\kappa \pi a^3}
e^{-\kappa \pi a^3} = 1$, $e^{\kappa \pi a^{3p}} e^{-\kappa \pi
a^{3p}} = 1$, $e^{\kappa \pi a^3} e^{-\kappa \pi a^{3p}} =
\theta(\kappa a^3) e^{\pi (\kappa^p - \kappa) a^{3p}}$, and
$e^{\kappa \pi a^3} e^{\kappa \pi a^{3p}} = \theta(\kappa a^3)
e^{\pi (\kappa^p + \kappa) a^{3p}}$, where $\theta$ is Dwork's first
splitting function on $\bb F_p$.

Using the compatibility relations, we may write
\[
V(a)^{-1} \f A_{a^2}^* V(a^p) a^{-(p-1)/2} =
\left(%
\begin{array}{cc}
  m_1 \theta(\kappa a^3) e^{\pi (\kappa^p - \kappa) a^{3p}} & m_2 \theta(\kappa a^3) e^{\pi (\kappa^p + \kappa) a^{3p}} \\
  m_2 \theta(-\kappa a^3) e^{-\pi (\kappa^p + \kappa) a^{3p}} & m_4 \theta(-\kappa a^3) e^{-\pi (\kappa^p - \kappa) a^{3p}} \\
\end{array}%
\right).
\]
Notice that the left-hand side belongs to $\hatR$, thus the
ride-hand side must be as well. It follows from Lemma \ref{L:
FermatLemma} below that:

\begin{theorem}
If $p \equiv 1,7 \> mod(12)$, then
\[
M=
\left(%
\begin{array}{cc}
  m_1 & 0 \\
  0 & m_4 \\
\end{array}%
\right),
\]
else if $p \equiv 5,11 \> mod(12)$, then
\[
M=
\left(%
\begin{array}{cc}
  0 & m_2 \\
  m_3 & 0 \\
\end{array}%
\right).
\]
\end{theorem}

\begin{lemma}\label{L: FermatLemma}
For $p \equiv 1,7 \> mod(12)$, we have $ord_p(\kappa^p -
\kappa)
> 0$ and $ord_p(\kappa^p + \kappa) = 0$. Otherwise, for $p
\equiv 5,11 \> mod(12)$, we have $ord_p(\kappa^p - \kappa) = 0$
and $ord_p(\kappa^p + \kappa) > 0$
\end{lemma}

\begin{proof}
Recall, $\kappa := \frac{2i}{3 \sqrt{3}}$. Thus,
\[
\kappa^p - \kappa =
\begin{cases}
i \left[ \left( \frac{2}{3 \sqrt{3}} \right)^p - \left( \frac{2}{3
\sqrt{3}} \right) \right] & \text{if } p \equiv 1 \text{ mod(4)} \\
-i \left[ \left( \frac{2}{3 \sqrt{3}} \right)^p + \left( \frac{2}{3
\sqrt{3}} \right) \right] & \text{if } p \equiv 3 \text{ mod(4)}.
\end{cases}
\]
Let us concentrate on $p \equiv 1$ mod(4). By Fermat's little
theorem, $2^p = 2 + p \c O_{\bb C_p}$ and $3^p = 3 + p \c O_{\bb
C_p}$, where $\c O_{\bb C_p} = \{ c \in \bb C_p : |c| \leq 1 \}$.
Consequently, we have
\[
\kappa^p - \kappa = \frac{6 \sqrt{3} i}{(3
\sqrt{3})^p}(1-3^{(p-1)/2}) + p \c O_{\bb C_p}.
\]
Thus, $\kappa^p-\kappa$ has positive $p$-adic order if and only if
$3^{(p-1)/2} \equiv 1$ mod($p$). By Euler's criterion, $3^{(p-1)/2}
\equiv \left(\frac{3}{p}\right)$ mod($p$), where the right-hand side
is the Legendre symbol. Since $\left(\frac{p}{3}\right) = 1$ if and
only if $p \equiv 1$ mod(3), the cases $p \equiv 1, 5$ mod(12) for
$\kappa^p - \kappa$ follow from quadratic reciprocity. The other
cases are similar.
\end{proof}

Notice that $v_1(-a) = v_2(a)$ and so $v_1(-a)' = -v_2(a)$. Thus, if
$T := \left(
\begin{array}{cc}
  0 & 1 \\
  1 & 0 \\
\end{array}
\right) = T^{-1}$, then $V(-a) T = V(a)$. Next, since the change $a
\mapsto -a$ does not affect $\f A_{a^2}$, changing $a \mapsto -a$ in
(\ref{E: DefFrob}) gives
\[
V(a)^{-1} \f A_{a^2} V(a^p) (-1)^{(p-1)/2} a^{-(p-1)/2} = S(a) T M T
S(a^p).
\]
Thus,
\[
(-1)^{(p-1)/2} TM = MT
\]
which demonstrates that if $p \equiv 1$ mod(12), then $m_1 = m_4$,
else if $p \equiv 7$ mod(12), then $-m_1 = m_4$. Similarly, if $p
\equiv 5$ mod(12), then $m_2=m_3$, else if $p \equiv 11$ mod(12)
then $-m_2 = m_3$.

Next, since $det(\f Y)$ is locally constant, from $\f A_{a^2} \f
Y^\phi = \f Y M$ we see that $det(M) = det(\f A_{a^2})$, and this
equals $p$ if $p \equiv 1$ mod(3) and $g := -g_2((p^2-1)/3)$ if $p
\equiv -1$ mod(3). Putting everything together, we have:
\[
M=
\left(%
\begin{array}{cc}
  \sqrt{p} & 0 \\
  0 & \pm \sqrt{p} \\
\end{array}%
\right) \quad \text{if } p \equiv 1,7 \text{ mod}(12)
\]
(of course, $\sqrt{p}$ may not be the positive square root) and
\[
M =
\left(%
\begin{array}{cc}
  0 & \sqrt{-g} \\
  \pm \sqrt{-g} & 0 \\
\end{array}%
\right) \quad \text{if } p \equiv 5,11 \text{ mod}(12)
\]
where ``$\pm$'' means positive if $p \equiv 1, 5 \> \text{ mod}(12)$
and negative if $p \equiv 7,11 \> \text{ mod}(12)$.

%%%-----------------------------------------------------------------------
\bigskip\noindent{\bf Symmetric powers at infinity.}
As demonstrated above, $\{y_1, y_2\}$ forms a formal basis for the
solution space (around a point near infinity) of the differential
operator $\tilde l := \frac{d^2}{da^2} - (1/a) \frac{d}{da} + 4\pi
a^4/3$. We define the $k$-th symmetric power of $\tilde l$, denoted
by $\tilde l_k$, as the unique monic differential equation of order
$k+1$ that annihilates every
\[
z_j := y_1^{k-j} y_2^j = a^{-k/2} e^{(k-2j) \kappa \pi a^3}
v_1(a)^{k-j} v_2(a)^j
\]
for $j = 0,1,\ldots, k$.

Define
\[
\hatR := \text{ analytic functions on } 1 < |a| < e, \text{ $e$
unspecified}.
\]
We are interested in studying the solutions of $\tilde l_k$ over the
ring $\hatR$. Let $Y(a) \in \hatR$ be a solution of $\tilde l_k$
which converges on $1 < |a| < e$. Fix $1 < |z| < e$. Since $\{ z_j
\}_{j=0}^k$ is a formal basis of $\tilde l_k$, if $Y$ is a solution
of $\tilde l_k$ over $\hatR$, then locally (in $D^-(z,1)$) using
Remark \ref{R: branch} we have the relation
\[
Y(a) = \sum_{j=0}^k c_j a^{-k/2} v^j(1/a) \bar v^{k-j}(1/a)[
e^{\kappa \pi (a-z)^3}]^{2j-k}
\]
for some constants $c_j \in \bb C_p$. Let us rewrite this as
\begin{equation}\label{E: poly}
a^{k/2} Y(a) = \sum_{r=0}^{p-1} S_r(a) [e^{\kappa \pi (a-z)^3}]^r
\end{equation}
where for each $r = 0 ,1, \ldots, p-1$,
\[
S_r(a) := \sum_{j \in E_r} c_j v^j(1/a) \bar v^{k-j}(1/a) [ e^{p
\kappa \pi (a-z)^3}]^{\frac{(2j-k)-r}{p}}
\]
with
\[
E_r := \{ j \in \bb Z \cap [0, k] : (2j-k) \equiv r \text{
mod}(p)\}.
\]
Since $v(1/a)$ and $\bar v(1/a)$ converge on $D^+(z,1) \subset
D^-(\infty, 1)$, and $e^{\pm p \kappa \pi (a-z)^3}$ and $Y(a)$
converge on $D^+(z,1)$, equation (\ref{E: poly}) shows us that the
functions $\exp(r \kappa \pi (a-z)^3)$, for $r=0, 1, \ldots, p-1$,
are linearly dependent over $\f M^+(z,1)[\sqrt{a}]$, where $\f
M^+(z,1)$ is the field of meromorphic functions on the closed unit
disc around $z$. However, Lemma \ref{L: irred} below shows us that
these functions are actually linearly independent over $\f
M^+(z,1)[\sqrt{a}]$. Consequently, we must have $S_r = 0$ for each
$r \geq 1$, and so,
\[
a^{k/2} Y(a) = S_0(a) = \sum_{(2j-k) \equiv 0 \text{ mod}(p)} c_j
v^j(1/a) \bar v^{k-j}(1/a) \exp[(2j-k) \kappa \pi (a-z)^3].
\]
Notice that the right-hand side now converges on $D^+(z,1)$. Thus,
if $k$ is odd, since $Y(a)$ converges on $D^+(z,1)$, we see that
$a^{k/2} \in \f M^+(z,1)$. However, by Lemma \ref{L: ExpSqrMero}
below, this cannot be. Therefore, for $k$ odd we must have $Y = 0$.

Thus, for $k$ even $\{ z_j : 0 \leq j \leq k, p | (k-2j) \}$ is a
basis for $ker(\tilde l_k | \hatR)$, else if $k$ is odd then
$ker(\tilde l_k | \hatR)=0$.

%%%-----------------------------------------------------------------------
\bigskip\noindent{\bf Symmetric powers of the Airy operator.} (cf.
\cite[\S 7]{pR86}) Define the Airy differential operator $l :=
\frac{d^2}{da^2} + \frac{\pi^2 a}{3}$ and denote by $l_k$ its $k$-th
symmetric power. We wish to study $ker(l_k | \hatR)$.

For $k$ odd, $ker(l_k | \hatR) = 0$ since any solution $Y \in \hatR$
of $l_k$ produces a solution $Y(a^2)$ to $\tilde l_k$. Assume $2 |
k$. Notice that $v_1(-a) = v_2(a)$. Consequently, for $0 \leq j \leq
k/2$ and $p | (k-2j)$,
\[
\widetilde{w}_j^+(a) := e^{(k-2j) \kappa \pi a^3} v_1(a)^{k-j}
v_2(a)^j + e^{-(k-2j)\kappa \pi a^3} v_1(a)^j v_2(a)^{k-j}
\]
is an even function while
\[
\widetilde{w}_j^-(a) := e^{(k-2j) \kappa \pi a^3} v_1(a)^{k-j}
v_2(a)^j - e^{-(k-2j)\kappa \pi a^3} v_1(a)^j v_2(a)^{k-j}
\]
is an odd function. Clearly, $\{ a^{-k/2} \widetilde{w}_j^+,
a^{-k/2} \widetilde{w}_j^- | 0 \leq j \leq k/2, p | (k-2j) \}$ is a
basis for $ker(\tilde l_k | \hatR)$.

From this, if $4 | k$ then $\{ a^{-k/2} \widetilde{w}_j^+(a) | 0
\leq j \leq k/2, p | (k-2j)\}$ is a set of even functions whereas
$\{ a^{-k/2} \widetilde{w}_j^-(a) | 0 \leq j \leq k/2, p | (k-2j)\}$
are odd. Thus $ker(l_k | \hatR) = span\{ w_j^+(a) | 0 \leq j \leq
k/2, p|(k-2j)\}$ where $w_j^+$ is defined by $w_j^+(a^2) = a^{-k/2}
\widetilde{w}_j^+(a)$. Similarly, if $2 | k$ but $4 \nmid k$, then
$ker(l_k | \hatR) = span\{ w_j^-(a) | 0 \leq j < k/2, p | (k-2j)
\}$, where $w_j^-$ is defined by $w_j^-(a^2) = a^{-k/2}
\widetilde{w}_j^-(a)$. Thus,
\[
dim_{\bb C_p} ker(l_k | \hatR) =
\begin{cases}
1 + \left\lfloor \frac{k}{2p} \right\rfloor & \text{if $4 | k$} \\
\left\lfloor \frac{k}{2p} \right\rfloor & \text{if $2 | k$ but $4
\nmid k$} \\
0 & \text{if $k$ odd}.
\end{cases}
\]

%%%-----------------------------------------------------------------------
\bigskip\noindent{\bf Appendix.}

\begin{lemma}\label{L: irred}
Let $\f M^+(0,1)$ denote the field of meromorphic functions on
$D^+(0,1)$ in the variable $a$. Then the polynomial $X^p - e^{p \pi
a} \in \f M^+(0,1)[X]$ is irreducible over $\f M^+(0,1)[\sqrt{a}]$.
\end{lemma}

\begin{proof}
It is clear that $X^p - e^{p \pi a}$ is either irreducible or
completely reducible over $\f M^+(0,1)[\sqrt{a}]$. If the latter,
then we would have
\[
e^{\pi a} = H_1(a) + \sqrt{a} H_2(a)
\]
for some $H_1, H_2 \in \f M^+(0,1)$. This implies that $e^{\pi a}$
satisfies a quadratic polynomial whose coefficients lie in $\c
M^+(0,1)$. Thus, $X^p-e^{p \pi a}$ is completely reducible over $\f
M^+(0,1)$. In particular, $e^{\pi a} \in \f M^+(0,1)$. However, this
is not true by the following lemma.
\end{proof}

\begin{lemma}\label{L: ExpSqrMero}
$\>$
\begin{itemize}
\item[(a)] $\exp(\pi a) \notin \f M^+(0,1)$.
\item[(b)] Let $a^{1/2}$ be a branch of square root centered around
$z \in \bb C_p^*$. Then $a^{1/2} \notin \f M^+(z, 1)$.
\end{itemize}
\end{lemma}

\begin{proof}
We know that $\exp(\pi a)$ converges on $D^-(0,1)$ and is divergent
on the boundary $|a| =1$. Suppose there exists analytic functions
$f(a)$ and $g(a)$ on $D^+(0,1)$ such that we have the equality
$e^{\pi a} = f(a)/g(a)$ whenever $|a| < 1$. Raising this to the
$p$-th power gives the equality $e^{p \pi a} = f(a)^p / g(a)^p$ for
all $|a| \leq 1$. Since $\exp(p \pi a)$ has no poles, $g(a)$ must be
invertible on $D^+(0,1)$. Let $h(a) := 1/g(a)$ and note this
converges on $D^+(0,1)$. Then $\exp(\pi a) = f(a) h(a)$ which shows
we may extend $\exp(\pi a)$ to an analytic function on $D^+(0,1)$.
Contradiction. This proves (a). A similar argument proves (b).
\end{proof}

%%% ----------------------------------------------------------------------
\section{Cohomological Formula for $M_k(T)$}\label{S: cohomForm}

Define
\[
M_k^*(T) := \prod_{\lambda \in |\bb G_m|} \prod_{i=0}^k (1 -
\pi_1(\lambda)^i \pi_2(\lambda)^{k-i} T^{deg(\lambda)})^{-1}
\]
where $|\bb G_m|$ denotes the set of Zariski closed points of
$\bb G_m := \bb P_{\bb F_p}^1 \setminus \{0, \infty\}$. In this
section, we will define a Frobenius operator $\bar \beta_k$ on
the (symmetric power) cohomology spaces $H_k^0$ and $H_k^1$
such that
\[
M_k^*(T) = \frac{det(1 - \bar \beta_k T | H_k^1)}{det(1-p \bar \beta_k
T | H_k^0)}.
\]
From this, we will deduce that $M_k(T)$ (and $M_k^*(T)$) is a
rational function. We will also find its field of definition.

The relationship between $M_k(T)$ and $M_k^*(T)$ is as follows.
If we write $M_k(T) = \exp \sum_{s\geq 1} N_s \frac{T^s}{s}$,
then taking the logarithmic derivative of the definition of
$M_k(T)$ shows
\[
N_s = \sum_{\bar a \in \bb F_{p^s}} \sum_{j=0}^k deg(\bar a)
[\pi_1(\bar a)^{k-j} \pi_2(\bar a)^j]^{\frac{s}{deg(\bar a)}}.
\]
Using the same procedure, we may calculate $N_s^*$, where
$M_k^*(T) = \exp \sum_{s\geq 1} N_s^* \frac{T^s}{s}$. It
follows that
\begin{equation}\label{E: Relations}
N_s = \sum_{j=0}^k (\pi_1(0)^{k-j} \pi_2(0)^j)^s + N_s^*.
\end{equation}

%%% ----------------------------------------------------------------------
\bigskip\noindent{\bf Differential Operator.} Define the operator
$\partial_a := a \frac{d}{da} + \pi a x$ on $\c K(b',b)$, and notice
that, formally
\[
\partial_a = e^{-\pi(x^3 + a x)} \circ a \frac{d}{da} \circ
e^{\pi(x^3 + a x)}.
\]
It follows that $\partial_a$ and $D_a$ commute, and so $\partial_a$
defines an operator on both $\c H_a(b',b)$ and $\c M_a(b',b)$.

Denote by $\Msk{a}(b',b)$ the $k$-th symmetric power of $\c
M_a(b',b)$ over $L(b')$. It is easy to see that $\partial_a$ induces
an endomorphism on $\Msk{a}(b',b)$ by extending linearly
\[
\partial_a(u_1 \cdots u_k) := \sum_{i=1}^k u_1 \cdots \hat u_i \cdots
u_k \partial_a(u_i)
\]
where $\hat u_i$ means we are leaving it out of the product.

%%%-----------------------------------------------------------------------
\bigskip\noindent{\bf Matrix representation of $\partial_a$.}
With the $L(b')$-basis $\{ v, w \}$ of $\c M_a(b',b)$, where $v :=
x$ and $w := x^2$, fix the ordered basis $(v^k, v^{k-1}w, \ldots,
w^k)$ of $\Msk{a}(b',b)$. Since $D_a(1) = 3 \pi x^3 + \pi a x$, we
see that $x^3 \equiv -\frac{ax}{3}$ in $\c M_a(b',b)$. Therefore,
$\partial_a(v) = \pi a w$ and $\partial_a(w) = \frac{-\pi a^2}{3}
v$. For each $j=0,1,\ldots, k$, we have
\begin{align*}
\partial_a( v^j w^{k-j}) &= j v^{j-1} \partial_a(v) w^{k-j} +
(k-j) v^j w^{k-j-1} \partial_a(w) \\
&= j \pi a v^{j-1}w^{k-j+1} + (k-j) \frac{-\pi a^2}{3}
v^{j+1}w^{k-j-1}.
\end{align*}
Thus, $\partial_a$, acting on \emph{row vectors}, has the matrix
representation $a \frac{d}{da} - G_k$ where
\[
G_k :=
\left(%
\begin{array}{ccccc}
  0 & k\pi a &  &  &  \\
  \frac{-\pi a^2}{3} & 0 & (k-1)\pi a &  &  \\
   & \frac{- 2 \pi a^2}{3} & & \ddots &  \\
   &  & \ddots &  & \pi a \\
   &  &  & \frac{-k \pi a^2}{3} & 0 \\
\end{array}%
\right).
\]

%%%-----------------------------------------------------------------------
\bigskip\noindent{\bf Dwork Operator.}
For $\frac{p-1}{p} \geq b > \frac{1}{p-1}$ and $b/p \geq b' > 0$,
$\bar \alpha(a)$ is an isomorphism of $\c M_a(b',b)$ onto $\c
M_{a^p}(b',b)$ satisfying $\partial_{a^p} \bar \alpha(a) = \bar
\alpha(a) \partial_a$ (since $a\frac{d}{da}$ and $\psi_x$ commute
with no $p$-factor). Thus $\bar \alpha(a)$ induces an isomorphism
$\Fsk{a}: \Msk{a}(b',b) \rightarrow \Msk{a^p}(b',b)$ defined by
acting on the image of $(u_1, \ldots, u_k) \in \c M_a(b',b)^k$ in
$\Msk{a}(b',b)$ by
\[
\bar \alpha_k(a)(u_1 \cdots u_k) := \bar \alpha(a)(u_1) \cdots \bar
\alpha(a)(u_k).
\]
Note, $\partial_{a^p} \Fsk{a} = \Fsk{a} \partial_a$.

Define the operator $\psi_a: L(b') \rightarrow L(p b')$ by
\[
\psi_a : \sum_{i=0}^\infty B_i a^i \longmapsto \sum_{i=0}^\infty
B_{pi} a^i.
\]
We may extend this map to $\Msk{a}(b',b)$ as follows. Fix the basis
$\{ v^i w^{k-i} \}_{i=0}^k$ of $\Msk{a}(b',b)$. Define the map
$\psi_a: \Msk{a^p}(b',b) \rightarrow \Msk{a}(p b',b)$ by
\[
\sum_{i=0}^k h_i v^i w^{k-i} \longmapsto \sum_{i=0}^k \psi_a(h_i)
v^i w^{k-i}.
\]
It is not hard to show that $\psi_a \partial_{a^p} = p
\partial_a \psi_a$.

Denote by $\Msk{a}(b)$ the space $\Msk{a}(b,b)$. We define the Dwork
operator $\beta_k := \psi_a \circ \Fsk{a}$ on $\Msk{a}(b)$ by noting
\[
\begin{CD}
\Msk{a}(b) @>\Fsk{a}>> \Msk{a}(b/p, b) @>\psi_a>> \Msk{a}(b).
\end{CD}
\]
This is a nuclear operator on $\Msk{a}(b)$ in the sense of
\cite{pR94}.

%%% ----------------------------------------------------------------------
\bigskip\noindent{\bf Cohomology.}
Define the cohomology spaces
\begin{align*}
H_k^0 &:= ker( \partial_a | \Msk{a}(b)) \\
H_k^1 &:= \Msk{a}(b) / \partial_a \Msk{a}(b).
\end{align*}
Since
\[
\beta_k \circ \partial_a = \psi_a \circ \Fsk{a} \circ
\partial_a = \psi_a \circ \partial_{a^p} \circ \Fsk{a} = p \partial_a \circ \psi_a \circ \Fsk{a} =
\partial_a \circ p \beta_k,
\]
$\beta_k$ induces endomorphisms
\[
\bar \beta_k : H^0_k \rightarrow H^0_k \quad \text{and} \quad \bar
\beta_k: H^1_k \rightarrow H^1_k.
\]
These maps are in fact isomorphisms; this follows by defining a
right-inverse to $\beta_k$ as in Remark \ref{R: nondeg}.

%%% ----------------------------------------------------------------------
\bigskip\noindent{\bf Dwork Trace Formula.}
From the following exact sequence and chain map:
\[
\begin{CD}
0 @>>> H_k^0 @>>> \Msk{a}(b) @>\partial_a>> \Msk{a}(b) @>>> H_k^1
@>>> 0 \\
@VVV @V p^s \bar\beta_k^s VV @V p^s \beta_k^s VV @VV \beta_k^s V @VV
\bar
\beta_k^s V @VVV \\
0 @>>> H_k^0 @>>> \Msk{a}(b) @>\partial_a>> \Msk{a}(b) @>>> H_k^1
@>>> 0
\end{CD}
\]
it follows from \cite[Proposition 8.3.10]{pR94} that the alternating
sum of the traces is zero:
\[
Tr(p^s \bar\beta_k^s | H_k^0) - Tr_{nuc}(p^s \beta_k^s | \Msk{a}(b))
+ Tr_{nuc}(\beta_k^s | \Msk{a}(b)) - Tr(\bar \beta_k^s | H_k^1) = 0.
\]
Hence,
\begin{equation}\label{E: HopfTrace}
(p^s-1) Tr_{nuc}(\beta_k^s | \Msk{a}(b)) = p^s Tr(\bar \beta_k^s |
H_k^0) - Tr(\bar \beta_k^s | H_k^1).
\end{equation}
Next, we have Dwork's trace formula
\begin{equation}\label{E: DworkTraceFormula}
(p^s - 1) Tr_{nuc}(\beta_k^s | \Msk{a}(b)) = \sum_{a \in \bb C_p^*,
a^{p^s-1}=1} Tr_{\Msk{a}(b)} \bar \alpha_k(a; s)
\end{equation}
where
\[
\bar \alpha_k(a;s) := \bar \alpha_k(a^{p^{s-1}}) \circ \bar
\alpha_k(a^{p^{s-2}}) \circ \cdots \circ \bar \alpha_k(a)
\]

%%% ----------------------------------------------------------------------
\bigskip\noindent{\bf Cohomological Formula for $M_k(T)$.}
Recall from \S \ref{S: Fibres}, if $\bar z \in \overline{\bb F}_p^*$
and $z$ is its Teichm\"uller representative in $\bb C_p$, then $\bar
\alpha_z$ was defined to be the operator $\bar \alpha(a)$
specialized at $a = z$. Using this, define
\[
\bar \alpha_{z,s}^{(k)} := \bar \alpha_k(a^{p^{s-1}}) \circ \bar
\alpha_k(a^{p^{s-2}}) \circ \cdots \circ \bar \alpha_k(a)|_{a=z}.
\]
Using (\ref{E: HopfTrace}) and (\ref{E: DworkTraceFormula}) for the
first two equalities (resp.) of the following, we have
\begin{align}\label{E: coh_form}
\frac{det(1 - \bar \beta_k T | H^1_k)}{det(1 - p \bar \beta_k T|
H^0_k)} &= \exp \sum_{s \geq 1} (p^s - 1) Tr_{nuc}(\beta_k^s |
\Msk{a}(b)) \frac{T^s}{s} \notag \\
&= \exp \sum_{s \geq 1} \sum_{a \in \bb C_p, a^{p^s-1}=1}
Tr_{\Msk{a}(b)}( \Fsk{a;s}) \frac{T^s}{s} \notag \\
&= \exp \sum_{r \geq 1} \sum_{\substack{\bar a \in \overline{\bb F}_p^*, deg(\bar a) = r \\
a = Teich(\bar a)}}
\sum_{m \geq 1} Tr_{\Msk{a}(b)}( \Fsk{a;rm}) \frac{T^{rm}}{rm} \notag \\
&= \exp \sum_{r \geq 1} \sum_{\substack{\bar a \in \overline{\bb F}_p^*, deg(\bar a) = r \\
a = Teich(\bar a)}}
\sum_{m \geq 1} Tr_{\bb C_p}( \bar \alpha_{a, rm}^{(k)}) \frac{T^{rm}}{rm} \notag \\
&= \exp \sum_{r \geq 1} \sum_{\substack{\bar a \in \overline{\bb F}_p^*, deg(\bar a) = r \\
a = Teich(\bar a)}}
\sum_{m \geq 1} Tr_{\bb C_p}( (\bar \alpha_{a, r}^{(k)})^m) \frac{T^{rm}}{rm} \notag \\
&= \prod_{r \geq 1} \prod_{\substack{\bar a \in \overline{\bb
F}_p^*, deg(\bar a) = r \\ a = Teich(\bar a)}} [det_{\bb C_p}(1-
\bar \alpha_{a,r}^{(k)} T^r )]^{-1/r}.
\end{align}
From \S \ref{S: Fibres} we know that $\bar \alpha_{a,r}$ has two
eigenvalues $\pi_1(a)$ and $\pi_2(a)$ on $M_a$. Since $\bar
\alpha_{a,r}^{(k)}$ is the $k$-th symmetric power of $\bar
\alpha_{a,r}$, its eigenvalues are $\pi_1(a)^{k-j} \pi_2(a)^j$ for
$j=0,1,\ldots, k$. Thus,
\[
det(1 - \bar \alpha_{a,r}^{(k)} T) = \prod_{j=0}^k (1 -
\pi_1(a)^{k-j} \pi_2(a)^j T).
\]
This means (\ref{E: coh_form}) equals
\[
\prod_{r \geq 1} \> \prod_{\substack{\bar a \in \overline{\bb F}_p^*, deg(\bar a) = r \\
a = Teich(\bar a)}} \> \prod_{j=0}^k (1 - \pi_1(a)^{k-j} \pi_2(a)^j
T^r)^{-1/r}
\]
which is precisely $M_k^*(T)$.

\begin{theorem}
$M_k(T)$ is a rational function with coefficients in $\bb
Z[\zeta_p]$. If $p \equiv 1 \> mod(3)$, then the coefficients of
$M_k(T)$ lie in the ring of integers of the unique subfield $L$ of
$\bb Q(\zeta_p)$ with $[L: \bb Q] = 3$. However, if $p \equiv -1 \>
mod(3)$ then $M_k(T)$ has integer coefficients.
\end{theorem}

\begin{proof}
We will first prove rationality by using a result of Borel-Dwork.
From the cohomological description, $M_k(T)$ is $p$-adic
meromorphic. Next, write $M_k(T) = \exp \sum_{s\geq 1} N_s
\frac{T^s}{s}$ as in the introduction, with
\[
N_s = \sum_{\bar a \in \bb F_{p^s}} \sum_{j=0}^k deg(\bar a)
[\pi_1(\bar a)^{k-j} \pi_2(\bar a)^j]^{\frac{s}{deg(\bar a)}}.
\]
Since $|\pi_j(\bar a)|_{\bb C} \leq p^{deg(\bar a)/2}$ by the
Riemann Hypothesis for curves, $|N_s|_{\bb C} \leq c s p^{(k+1)
s/2}$ for a positive integer $c$ independent of $s$. Thus $\sum_{s
\geq 1} \frac{N_s}{s} T^s$ has at least a $p^{-(k+1)/2}$ radius of
convergence on $\bb C$. Consequently, $M_k(T)$ has a positive $\bb
C$-radius of convergence.

Let $\zeta_p$ is a primitive $p$-th root of unity. Since $L(x^3 +
\bar a x, T)$ is a polynomial with entries in $\bb Z[\zeta_p]$, we
see that the roots live in $\bb Z[\zeta_p]$-rational orbits under
the action of the Galois group $Gal(\overline{\bb Q}/\bb
Q(\zeta_p))$. Hence, the symmetric polynomial $\prod_{j=0}^k (1 -
\pi_1(\bar a)^{k-j} \pi_2(\bar a)^j T)$ must also have coefficients
in $\bb Z[\zeta_p]$. It follows that $M_k(T)$ has coefficients in
$\bb Z[\zeta_p]$. Since $M_k(T)$ is $p$-adic meromorphic and
converges on a disc of positive radius in $\bb C$, and the
coefficients lie in a fixed number field, by the Borel-Dwork theorem
\cite[Theorem 3]{bD60} $M_k(T)$ is a rational function with
coefficients in $\bb Z[\zeta_p]$.

Next, with each $\lambda \in \bb F_p^*$ we may define
$\sigma_\lambda \in Gal(\bb Q(\zeta_p) / \bb Q)$, where , by
$\sigma_\lambda(\zeta_p) := \zeta_p^\lambda$; note, every element of
$Gal(\bb Q(\zeta_p) / \bb Q)$ may be represented uniquely by such a
$\sigma_\lambda$. Now,
\begin{align*}
\sum_{\bar x \in \bb F_{p^k}} \zeta_p^{Tr_{\bb F_{p^k} / \bb
F_p}(\bar x^3 + \lambda^2 \bar a\bar x)} &= \sum_{\bar x \in \bb
F_{p^k}} \zeta_p^{Tr_{\bb
F_{p^k}/ \bb F_p}[\lambda^3(\bar x^3 + \bar a\bar x)]} \\
&= \sigma_{\lambda^3} \left( \sum_{\bar x \in \bb F_{p^k}}
\zeta_p^{Tr_{\bb F_{p^k} / \bb F_p}(\bar x^3 + \bar a\bar x)}
\right).
\end{align*}
This means $\sigma_{\lambda^3}(L(x^3+\bar a x, T)) = L(x^3 +
\lambda^2 \bar a x, T)$. Consequently,
\[
\sigma_{\lambda^3}(M_k(T)) = M_k(T),
\]
and so $M_k(T)$ is fixed by the subgroup $\{ \lambda^3 | \lambda \in
\bb F_p^*\}$. Since this subgroup has index $r := gcd(3,p-1)$ in the
Galois group, the coefficients of $M_k(T)$ must lie in the ring of
integers of the unique subfield of $\bb Q(\zeta_p)$ of degree $r$
over $\bb Q$.
\end{proof}

%%% ----------------------------------------------------------------------
\section{Cohomology of Symmetric Powers}\label{S: CohomSymPow}

%%%-----------------------------------------------------------------------
\subsection{$H_k^0$ for all odd $k$, or for all $k < 2p$ even}

Recall, with respect to the basis $\{ x, x^2\}$ of $\c M_a(b',b)$,
$\partial_a$ has the matrix representation (acting on \emph{row}
vectors):
\[
\partial_a = a \frac{d}{da} -
\left(%
\begin{array}{cc}
  0 & \pi a \\
  \frac{-\pi a^2}{3} & 0 \\
\end{array}%
\right).
\]
Thus, if $(C_1, C_2)$ is a solution of this differential equation
then $C_2$ must satisfy the Airy differential equation: $y'' +
\frac{\pi^2 a}{3} y = 0$.

\begin{theorem}\label{T: H_0}
Let $k$ be a positive integer. If $k$ is odd, or $k$ is even with $k
< 2p$, then $H_k^0 := ker(\partial_a|\Msk{a}(b',b)) = 0$.
\end{theorem}

\begin{proof}
Every element of $\Msk{a}(b',b)$ which is annihilated by
$\partial_a$ corresponds to a unique solution in $L(b')$, and hence
$\hatR$, of the scalar equation $l_k$, where $l_k$ was defined in \S
\ref{S: AiryAsymp}. By the methods of \S \ref{S: AiryAsymp}, if $k$
is odd, then no such solution exists in $\hatR$. If $2 | k$ but $4
\nmid k$, and $k < 2p$, then again, no such solution exists. If $4 |
k$, then one solution of $l_k$ exists in $\hatR$, however, this
corresponds to $v_1(a)^{k/2} v_2(a)^{k/2}$ which is clearly not an
element of $L(b')$.
\end{proof}

\begin{conjecture}
$H_k^0 = 0$ for all positive integers $k$.
\end{conjecture}

%%% ----------------------------------------------------------------------
\subsection{$H_k^1$ for odd $k < p$}\label{S: oddCohom}

We will assume $k$ is odd and $k < p$ throughout this section. The
two vectors $v := x$ and $w := x^2$ form an $L(b)$-basis of $\c
M_a(b)$. Thus, the set $\{ v^j w^{k-j} | 0 \leq j \leq k \}$ is an
$L(b)$-basis of $\Msk{a}(b)$. With $\delta_j := \frac{bj}{3}$,
define
\begin{align*}
\Msk{a}(b; \bs \delta + \rho) &:= \bigoplus_{j=0}^{k} L(b;
\delta_{j+1} + \rho) v^{k-j}w^{j} \\
&= L(b; \delta_1 + \rho) \oplus L(b; \delta_2 + \rho) \oplus \cdots
\oplus L(b; \delta_{k+1} + \rho)
\end{align*}
and
\[
\Msk{a}(b) := \bigcup_{\rho \in \bb R} \Msk{a}(b; \bs \delta +
\rho).
\]
We will write $\Msk{a}(b;\rho)$ with the $\bs \delta$ implied.
Define
\begin{align}\label{E: DecompV}
V_k &:= \bb C_p^{k+1} \oplus \bigoplus_{j=0}^{(k-1)/2} (\bb C_p a) v^{k-2j} w^{2j} \\
V_k(b;\rho) &:= V_k \cap \Msk{a}(b;\rho) \notag \\
V_k(b) &:= \bigcup_{\rho \in \bb R} V_k(b; \rho). \notag
\end{align}

\begin{theorem}\label{T: oddDecomp}
Let $p \geq 5$ be a prime. Suppose $k$ is odd and $k < p$. Let $b$
be a real number such that $e := b - \frac{1}{p-1} > 0$. Then we
have the following decomposition:
\[
\Msk{a}(b) = V_k(b) \oplus \partial_a \Msk{a}(b).
\]
Furthermore, $\partial_a$ is injective on $\Msk{a}(b)$.
\end{theorem}

It will be useful to recall, if we order the basis of $\Msk{a}(b)$
as $(v^k, v^{k-1}w, \ldots, w^k)$, then $\partial_a$ has the matrix
representation $\partial_a = a \frac{d}{da} - G_k$, acting on
\emph{row vectors}, where
\[
G_k :=
\left(%
\begin{array}{ccccc}
  0 & k\pi a &  &  &  \\
  \frac{-\pi a^2}{3} & 0 & (k-1)\pi a &  &  \\
   & \frac{- 2 \pi a^2}{3} & & \ddots &  \\
   &  & \ddots &  & \pi a \\
   &  &  & \frac{-k \pi a^2}{3} & 0 \\
\end{array}%
\right)
\]
That is, if we write the row vector $\xi = (\xi_1, \ldots,
\xi_{k+1})$ for $\sum_{j=1}^{k=1} \xi_j v^{k-j} w^j$ then
\[
\partial_a \xi \quad \text{means} \quad (\xi_1, \ldots, \xi_{k+1})(a
\frac{d}{da} - G_k).
\]

\begin{lemma}
$\Msk{a}(b;0) \subset V_k(b;0) + \Msk{a}(b;e) G_k$.
\end{lemma}

\begin{proof}
{\bf Odd Column Case.} For $n \geq 2$, consider the vector $(0,
\ldots, B_n a^n, \ldots, 0)$ where $a^n$ is in the $(2j+1)$-entry
with all others zero, and $B_n a^n \in L(b; \delta_{2j+1})$.

Suppose $n \leq (k+1)/2 - j$. Then
\[
(0, \ldots, \overbrace{ B_n a^n}^{2j+1 \text{ entry}}, \ldots, 0) =
(0, 0, \ldots, 0, C_{2(j+1)}, 0, C_{2(j+2)}, \cdots) G_k +
(0,\ldots, \overbrace{\eta_{j,n}}^{2(n-1+j)+1 \text{ entry}},
\ldots, 0)
\]
where
\[
\eta_{j,n} := \prod_{l=0}^{n-2} \frac{3 [k-(2(j+l)+1)]}{[2(j+l)+1]}
B_n a \in L(b; \delta_{2(n-1+j)+1})
\]
and for $i=1, 2, \ldots, n-1$, we have
\[
C_{2(j+i)} = \frac{-3^i}{\pi} \left( \frac{\prod_{l=0}^{i-2}
[k-(2(j+l)+1)]}{\prod_{l=0}^{i-1} [2(j+l)+1]} \right) B_n a^{n-1-i}
\in L(b; \delta_{2(j+i)} + e)
\]

Next, suppose $n \geq (k+1)/2-j+1$. Then
\[
(0, \ldots, \overbrace{ B_n a^n}^{2j+1 \text{ entry}}, \ldots, 0) =
(0, 0, \ldots, 0, C_{2(j+1)}, 0, C_{2(j+2)}, \cdots) G_k
\]
where, for $m = 1, 2 \ldots, (k+1)/2 - j$,
\[
C_{2(j+m)} := \frac{(-3)^m
(j-\frac{k-1}{2})_{m-1}}{2\pi(j+\frac{1}{2})_m} B_n a^{n-m-1}.
\]
Notice that, since $k < p$, the constants in front of the series are
$p$-adic units.

{\bf Even Column Case.} Fix $j \in \{1, \ldots, \frac{k+1}{2} \}$.
Let $g(a) := \sum_{i=0}^\infty B_i a^i \in L(b; \delta_{2j})$. Then
we may write
\[
(0, \ldots, 0, \overbrace{g(a)}^{\text{entry }2j}, 0, \ldots, 0) =
\xi G_k + (0, \ldots, 0, \overbrace{B_0}^{\text{entry }2j}, 0,
\ldots, 0)
\]
where
\[
\xi := (\xi_1, 0, \xi_3, 0, \ldots, \xi_{2j-1}, 0, 0, \ldots, 0) \\
\]
where, for $m=0,1, \ldots, j-1$,
\[
\xi_{2j-1-2m} := \frac{(1-j)_m}{2(-3)^m \pi (\frac{k}{2}-j+1)_{m+1}}
\sum_{i=0}^\infty B_{i+1} a^{i+m} \in L(b; \delta_{2j-1-2m} + e).
\]
\end{proof}

\begin{lemma}
$V_k(b) \cap  \left(\Msk{a}(b) G_k \right) = \{0\}$.
\end{lemma}

\begin{proof}
We will proceed by induction. Let $\eta = (\eta_1, \ldots,
\eta_{k+1}) \in V_k \cap \left(\Msk{a}(b) G_k \right)$. Then there
is a $\xi \in \Msk{a}(b)$ such that $\xi G_k = \eta$. In particular,
$\eta_1 = \frac{-\pi a^2}{3} \xi_2$. Hence, $\xi_2 = 0$. Suppose
$\xi_{2(j-1)} = 0$, then for some $p$-adic units $c_1$ and $c_2$ we
may write
\[
\eta_{2j-1} = c_1 a \xi_{2(j-1)} + c_2 a^2 \xi_{2j}
\]
which shows $\xi_{2j} = 0$. Thus, $\xi_{2j} = 0$ for $j=1, 2,
\ldots, \frac{k+1}{2}$. Next, since $\eta_{k+1} = \pi a \xi_k$ we
see that $\xi_k = 0$. Suppose $\xi_{2j+1} = 0$, then since there are
units $c_1$ and $c_2$ such that
\[
\eta_{2j} = c_1 a \xi_{2j-1} + c_2 a^2 \xi_{2j+1}
\]
we see that $\xi_{2j-1} = 0$.
\end{proof}

\begin{lemma}\label{L: Glift}
Let $u \in \Msk{a}(b)$ such that $u G_k \in \Msk{a}(b;\rho)$. Then
$u \in \Msk{a}(b; \rho + e)$. Consequently, $G_k$ is an injective
operator on $\Msk{a}(b)$.
\end{lemma}

\begin{proof}
We will proceed by induction. Let $\xi := u G_k$. Then $\xi_1 =
\frac{-\pi}{3} a^2 u_2 \in L(b; \delta_1)$. Thus, $u_2 \in L(b;
\delta_2+ e)$ as desired. Next, suppose $u_{2j} \in L(b; \delta_{2j}
+ e)$. Then $\pi a u_{2j} \in L(b; \delta_{2j+1})$. From $\xi = u
G_k$, there are $p$-adic units $c_1$ and $c_2$ such that
\[
\xi_{2j+1} = c_1 \pi a u_{2j} - \frac{c_2 \pi a^2}{3} u_{2j+2} \in
L(b; \delta_{2j+1}).
\]
Hence, $\frac{c_2 \pi a^2}{3} u_{2j+2} \in L(b; \delta_{2j+1})$, and
this implies $u_{2j+2} \in L(b; \delta_{2j+2} + e)$.

Next, by hypothesis, $\xi_{k+1} = \pi a u_k \in L(b; \delta_{k+1})$
and so, $u_k \in L(b; \delta_k + e)$. Suppose $u_{2j+1} \in L(b;
\delta_{2j+1} + e)$. Then $\pi a^2 u_{2j+1} \in L(b;\delta_{2j})$.
As before, from $\xi = u G_k$, there are $p$-adic units $c_1$ and
$c_2$ such that
\[
\xi_{2j} = c_1 \pi a u_{2j-1} - c_2 \pi a^2 u_{2j+1} \in L(b;
\delta_{2j}).
\]
As before, it follows that $u_{2j-1} \in L(b; \delta_{2j-1} + e)$.
\end{proof}

\begin{lemma}\label{L: SymDecomp_D}
$\Msk{a}(b; 0) \subset V_k(b;0) + \partial_a \Msk{a}(b; e)$.
\end{lemma}

\begin{proof}
Let $u \in \Msk{a}(b;0)$. Set $u^{(0)} := u$. We know there exists
unique $\eta^{(0)} \in V_k(b;0)$ and $\xi^{(0)} \in \Msk{a}(b; e)$
such that
\[
u^{(0)} = \eta^{(0)} + \xi^{(0)} G_k.
\]
Set $E := a \frac{d}{da}$. Then
\begin{align*}
u^{(1)} :&= u^{(0)} - \eta^{(0)} + \partial_a \xi^{(0)} \\
&= E \xi^{(0)} \in \Msk{a}(b;e).
\end{align*}
Again, there exists unique $\eta^{(1)} \in V_k(b;e)$ and $\xi^{(1)}
\in \Msk{a}(b;2e)$ such that
\[
u^{(1)} = \eta^{(1)} + \xi^{(1)} G_k.
\]
Set $u^{(2)} := u^{(1)} - \eta^{(1)} + \partial_a \xi^{(1)}$.
Continue this procedure $h$-times:
\[
u^{(h)} := u^{(h-1)} - \eta^{(h-1)} + \partial_a \xi^{(h-1)} \in
\Msk{a}(b;he).
\]
Adding these together produces
\[
u^{(h)} = u^{(0)} - \sum_{i=0}^{h-1} \eta^{(i)} + \partial_a
\sum_{i=0}^{h-1} \xi^{(i)}.
\]
Letting $h \rightarrow \infty$ we see that $u^{(h)} \rightarrow 0$
in $\Msk{a}(b)$. Thus,
\[
u^{(0)} = \sum_{i=0}^\infty \eta^{(i)} - \partial_a
\sum_{i=0}^\infty \xi^{(i)} \in V_k(b) + \partial_a \Msk{a}(b).
\]
\end{proof}

\begin{lemma}
Let $u \in \Msk{a}(b)$ be such that $\partial_a u \in
\Msk{a}(b;\rho)$. Then $u \in \Msk{a}(b; \rho + e)$.
\end{lemma}

\begin{proof}
Suppose $u \ne 0$. Choose $c \in \bb R$ such that $u \in
\Msk{a}(b;c)$ but $u \notin \Msk{a}(b;c+e)$. By hypothesis,
$\partial_a u \in \Msk{a}(b;\rho)$. Thus
\[
- u G_k = \partial_a u - E u \in \Msk{a}(b; \rho) + \Msk{a}(b;c) =
\Msk{a}(b; l).
\]
where $l := \min\{\rho, c\}$. Hence, $u \in \Msk{a}(b;l + e)$ which
means $l = \rho$.
\end{proof}

\begin{corollary}
$\partial_a$ is an injective operator on $\Msk{a}(b)$.
\end{corollary}

\begin{lemma}
$V_k(b) \cap \partial_a \Msk{a}(b) = \{ 0 \}$.
\end{lemma}

\begin{proof}
Let $u \in V_k(b) \cap \partial_a \Msk{a}(b)$. Assume $u \ne 0$.
Then, we may find $\rho \in \bb R$ such that $u \in \Msk{a}(b;\rho)$
but $u \notin \Msk{a}(b; \rho+e)$. By hypothesis, there exists $\eta
\in \Msk{a}(b)$ such that $\partial_a \eta = u$. Thus, $\eta \in
\Msk{a}(b;\rho+e)$. Now,
\[
u + \eta G_k = E \eta \in \Msk{a}(b;\rho+e)
\]
which means that we may find $\zeta \in V_k(b;\rho+e)$ and $w \in
\Msk{a}(b; \rho+2e)$ such that
\[
E \eta = \zeta + w G_k.
\]
Therefore, $u = \zeta + (w-\eta)G_k$. Since $u \in V_k(b)$, by
uniqueness, $u = \zeta$ (and $\eta = w$). However, this means $u \in
V_k(b; \rho+e) \subset \Msk{a}(b; \rho+e)$ which is a contradiction.
\end{proof}

This concludes the proof of Theorem \ref{T: oddDecomp}.

%%% ----------------------------------------------------------------------
\subsection{$H_k^1$ for even $k < p$}\label{S: EvenDecomp}

% ----------------------------------------------------------------
\newcommand{\Mk}{\mathcal M^{(k)}_a}
\newcommand{\tilMk}{\mathcal{\widetilde{M}}^{(k)}_a}
% ----------------------------------------------------------------

In this section, while we attempt to proceed in a similar manner to
the last, we are only able to prove a partial decomposition theorem.
The crucial difference is that when $k$ is even, $G_k$ has a kernel.

Assume $k$ is even and $k < p$. With $\delta_j := \frac{bj}{3}$,
define
\[
\Mk(b;\bs \delta + \rho) := \bigoplus_{j=0}^{k}L(b; \rho +
\delta_{j+1}) v^{k-j} w^{j}
\]
As we did in the previous section, we will write $\Msk{a}(b;\rho)$
with the $\bs \delta$ implied. Ordering the basis as $(v^k, v^{k-1}
w, \ldots, w^k)$, define the subspace $\tilMk(b; \rho) \subset
\Mk(b;\rho)$ as
\[
\tilMk(b;0) := L(b;\delta_1) \oplus L(b; \delta_2) \oplus \cdots
\oplus L(b; \delta_k) \oplus \{0\}
\]
Next, define
\[
V_k := \bb C_p^{k+1} \oplus \bigoplus_{j=0}^{k/2-1} \bb C_p a
v^{k-2j} w^{2j},
\]
and the vector
\[
{\bf k} := (K_1, 0, K_3, 0, \ldots, 0, K_{k+1})
\]
where for $j=0,1, \ldots, k/2$
\[
K_{2j+1} := \frac{\prod_{i=1}^j [\frac{k}{2} - (i-1)]}{j! \>
3^{k/2-j}} a^{k/2-j}.
\]
Notice that the kernel of $G_k$ is $L(b) {\bf k}$. Define
\[
\c I_{ker}(b; \rho) = a L(b) {\bf k} \cap \Mk(b; \rho) \quad
\text{and} \quad \c I_{ker}(b) := a L(b) {\bf k}.
\]

\begin{theorem}\label{T: evenDecompThm}
Let $p \geq 5$. Suppose $k$ is even and $k < p$. Let $b$ be a real
number such that $b - \frac{1}{p-1} > 0$. Then
\[
\Msk{a}(b) = V_k(b) \oplus \partial_a \tilMk(b) \oplus \c
I_{ker}(b).
\]
\end{theorem}

The proof will consist of several lemmas.

\begin{lemma}\label{T: evenDecomp}
Set $e := b - \frac{1}{p-1} > 0$. Then
\[
\Mk(b;0) \subset V_k(b;0) + \tilMk(b;e) G_k + \c I_{ker}(b;0).
\]
\end{lemma}

The proof of this is a bit involved.

\bigskip\noindent{\bf Even Entry Case.}
Consider $(0, \ldots, g(a), \ldots, 0)$ where $g(a) = \sum_{i \geq
0} B_i a^i \in L(b; \delta_{2j})$ is in the $2j$-entry, and all
other entries are zero. Then
\[
(0, \ldots, \overbrace{g(a)}^{2j \text{ entry}}, \ldots, 0) = (C_1,
0, C_3, 0, \ldots, C_{2j-1}, 0, 0, \ldots, 0) G_k + (0, \ldots,
\overbrace{B_0}^{2j \text{ entry}}, \ldots, 0)
\]
where
\[
C_{2j-2l-1} = \frac{ \prod_{m=0}^{l-1} (j-l+m)}{2 \pi 3^l
\prod_{m=1}^{l+1} (\frac{k}{2} - j + m)} \sum_{i \geq 1} B_i
a^{i-1+l} \in L(b; \delta_{2j-2l-1} + e)
\]
for $l = 0, 1, \ldots, j-1$ and all other coordinates are zero.
Note, $C_{k+1}$ equals zero.

\bigskip\noindent{\bf Odd Entry Case.} This is a bit more technical.
For $n \geq 2$, consider the vector $(0, \ldots, B_n a^n, \ldots,
0)$ where $a^n$ is in the $(2j+1)$-entry with all others zero, and
$B_n a^n \in L(b; \delta_{2j+1})$.

First, suppose $n \leq k/2 - j$. Then
\[
(0, \ldots, \overbrace{ B_n a^n}^{2j+1 \text{ entry}}, \ldots, 0) =
(0, 0, \ldots, 0, C_{2(j+1)}, 0, C_{2(j+2)}, \cdots) G_k +
(0,\ldots, \overbrace{\eta_{j,n}}^{2(n-1+j)+1 \text{ entry}},
\ldots, 0)
\]
where
\[
\eta_{j,n} := \prod_{l=0}^{n-2} \frac{3 [k-(2(j+l)+1)]}{[2(j+l)+1]}
B_n a \in L(b; \delta_{2(n-1+j)+1})
\]
and for $i=1, 2, \ldots, n-1$, we have
\[
C_{2(j+i)} = \frac{-3^i}{\pi} \left( \frac{\prod_{l=0}^{i-2}
[k-(2(j+l)+1)]}{\prod_{l=0}^{i-1} [2(j+l)+1]} \right) B_n a^{n-1-i}
\in L(b; \delta_{2(j+i)} + e)
\]

Next, suppose $n \geq \frac{k}{2} - j+1$. We will demonstrate that
\begin{equation}\label{E: oddReduct}
(0, \ldots, \overbrace{B_n a^n}^{2j+1 \text{ entry}}, \ldots, 0) =
(0,C_2,0,C_4,0, \ldots, C_k,0) G_k + {\bf k} h
\end{equation}
where $h$ is given by (\ref{E: h}) below and, for $l = 1, 2, \ldots,
k/2$,
\[
C_{2l} = \frac{\bb Z}{2^k (k/2)! \pi} B_n a^{n+j-l-1} \in L(b;
\delta_{2j} + e)
\]
where ``$\bb Z$'' indicates some determinable integer.

In (\ref{E: oddReduct}), since the even numbered entries are all
zero, let us \emph{delete them} from this equation to obtain
\begin{equation}\label{E: oddReduct2}
(0, \ldots, \overbrace{B_n a^n}^{j^\text{th} \text{ entry}}, \ldots,
0) = (C_2, C_4, \ldots, C_k) \widetilde{G}_k + {\bf k}_{odd} h
\end{equation}
where the vector on the left has $(k/2)+1$ entries and $B_n a^n$ is
in the $j$-th entry, ${\bf k}_{odd}$ is the vector ${\bf k}$ with
the even entries deleted, and $\widetilde{G}_k$ is the
$(k/2)\times(k/2+1)$-matrix
\[
\widetilde{G}_k = \left(
  \begin{array}{ccccc}
    -\frac{\pi a^2}{3} & (k-1)\pi a &  &  &  \\
     & -\frac{3\pi a^2}{3} & (k-3) \pi a &  &  \\
     &  & \ddots & \ddots &  \\
     &  &  & -\frac{(k-1)\pi a^2}{3} & \pi a \\
  \end{array}
\right).
\]
Thus, we may rewrite this as
\[
(0, \ldots, \overbrace{B_n a^n}^{j+1 \text{ entry}}, \ldots, 0) =
(h, C_2, C_4, \ldots, C_k) N_k \quad \text{where} \quad N_k :=
\left(
  \begin{array}{c}
    {\bf k}_{odd} \\
    \widetilde{G}_k \\
  \end{array}
\right).
\]
Note, $N_k$ is a square matrix. Furthermore, $det N_k = \frac{2^k
(k/2)!}{3^{k/2}} \pi^{k/2} a^k$ by the following lemma.

\begin{lemma}
$det N_k = \frac{2^k (k/2)!}{3^{k/2}} \pi^{k/2} a^k$
\end{lemma}

\begin{proof}
Expanding the determinant always using the left column, we get
\[
det N_k = \frac{\pi^{k/2} a^k}{3^{k/2}} \sum_{j=0}^{k/2}
\frac{k!}{2^{k/2} \cdot (k/2)!} \binom{k/2}{j}^2 \binom{k}{2j}^{-1}.
\]
Thus, the lemma follows from the following general formula:
\begin{equation}\label{E: Combo}
\sum_{j=0}^n \binom{n}{j}^2 \binom{2n}{2j}^{-1} = \frac{(2^n \cdot
n!)^2}{(2n)!}.
\end{equation}
To prove this last formula (of which I thank Zhi-Wei Sun for the
following), notice that the above may be rewritten as
\[
\sum_{j=0}^n \binom{n}{j}^2 \binom{2n}{2j}^{-1} = 4^n
\binom{2n}{n}^{-1}.
\]
Therefore, since
\[
\binom{2n}{n} \binom{n}{j}^2 \binom{2n}{2j}^{-1} =
\binom{2j}{j}\binom{2n-2j}{n-j}
\]
we have (\ref{E: Combo}) equal to
\[
\sum_{j=0}^n \binom{2j}{j} \binom{2n-2j}{n-j}=4^n
\]
which is a well-known combinatorial formula.
\end{proof}

By Cramer's rule, with $N_k(1)$ denoting the matrix $N_k$ with the
first row replaced by $(0, \ldots, B_n a^n, \ldots, 0)$, where $B_n
a^n$ is in the $(j+1)$-entry, we have
\begin{equation}\label{E: h}
h = \frac{det N_k(1)}{det N_k} = \frac{3^{\frac{k}{2}-j}}{2^k
(k/2)!} \left( \prod_{i=1}^j (2i-1) \right) \left(
\prod_{i=j}^{\frac{k}{2}-1} [ k-(2i+1)] \right) B_n a^{n -
\frac{k}{2}+j}.
\end{equation}
Notice that $h$ is divisible by $a$ by our hypothesis that $n \geq
\frac{k}{2} - j+1$. Using this and (\ref{E: oddReduct2}), it follows
that ${\bf k} h \in \c I_{ker}(b;0)$ and
\[
C_{2l} = \frac{\bb Z}{2^k (k/2)! \pi} a^{n+j-l-1} \qquad \text{for }
l = 0, 1, 2, \ldots, k/2.
\]
This finishes the proof of Lemma \ref{T: evenDecomp}.

\begin{lemma}
We have
\begin{enumerate}
\item  $V_k \cap \left( \Mk(b) G_k \right) = \{0\}$
\item $\c I_{ker}(b) \cap V_k = \{ 0 \}$
\item $\c I_{ker}(b) \cap \left( \tilMk(b) G_k \right) = \{ 0 \}$
\end{enumerate}
\end{lemma}

\begin{proof}
To prove the first statement, consider $\eta \in V_k \cap \left(
\Mk(b) G_k \right)$. Then $\eta_1 = \frac{-\pi a^2}{3} C_2$. Thus,
$\eta_1$ is divisible by $a^2$ which means $\eta_1 = C_2 = 0$. Next,
$\eta_3 = (k-1) \pi a C_2 - \pi a^2 C_4$. Thus, $C_4=0$. Continuing,
we see that $C_{2j} = 0$ for each $j$. Consequently, all the odd
coordinates of $\eta$ are zero. Next, since all the even coordinates
of $\eta$ are constants and the image of $\Mk(b)$ by $G_k$ is
divisible by $a$, $\eta$ must be zero.

The second statement is clear.

Lastly, let us prove the third. Suppose there exists $h \in L(b)$
such that $(C_1, C_2, \ldots, C_k, 0) G_k = h {\bf k}$. Looking at
the last entry in this, we must have $C_k = \frac{1}{\pi a} (\bb
Q_{>0}) h $ where ``$\bb Q_{>0}$'' is some positive rational number.
Working backwards from this using entries $k-1, k-3, \ldots, 3$, we
see $C_{k-2l} = \frac{1}{\pi a} ( \bb Q_{>0} ) h a^l$ for $l=0, 1,
\ldots, k/2-1$. However, the first entry says $\frac{-\pi a^2}{3}
C_2 = h \frac{1}{3^{k/2}} a^{k/2}$ which means $C_2 = -\frac{1}{\pi
a} (\bb Q_{> 0}) h a^{k/2}$. This contradicts the positivity of the
rational constant.
\end{proof}

\begin{corollary}\label{C: evenDirect}
$\Mk(b) = V_k \oplus \tilMk(b) G_k \oplus \c I_{ker}(b).$
\end{corollary}

\begin{lemma}
If $u \in \tilMk(b)$ and $u G_k \in \Mk(b;\rho)$, then $u \in
\tilMk(b;\rho+e)$.
\end{lemma}

\begin{proof}
We will proceed by induction. Let $\xi = u G_k$. Then $\xi_1 =
-\frac{\pi a^2}{3} u_2 \in L(b; \delta_1)$. Thus, $u_2 \in L(b;
\delta_2 + e)$. Next, suppose $u_{2j} \in L(b; \delta_{2j} + e)$.
Then there are $p$-adic units $c_1$ and $c_2$ such that
\[
\xi_{2j+1} = c_1 \pi a u_{2j} + c_2 \pi a^2 u_{2j+2} \in L(b;
\delta_{2j+1}).
\]
Since $c_1 \pi a u_{2j} \in L(b; \delta_{2j+1})$, $u_{2j+2} \in L(b;
\delta_{2j+2} + e)$ as desired. This finishes the \emph{odd}
coordinates of $\xi$. We now move on to the even.

By the form of $G_k$, we can find $p$-adic units $c_1$ and $c_2$
such that
\[
\xi_k =  c_1 \pi a u_{k-1} + c_2 \pi a^2 u_{k+1} \in L(b;
\delta_{k}).
\]
However, since $u \in \tilMk(b)$, we have $u_{k+1} = 0$. Hence,
$u_{k-1} \in L(b; \delta_{k-1} + e)$. Next, suppose $u_{2j+1} \in
L(b; \delta_{2j+1}+e)$. Again, we may find $p$-adic units $c_1$ and
$c_2$ such that
\[
\xi_{2j} = c_1 \pi a u_{2j-1} + c_2 \pi a^2 u_{2j+1} \in L(b;
\delta_{2j}).
\]
It follows that $u_{2j-1} \in L(b; \delta_{2j-1}+e)$.
\end{proof}

Proceeding as in the previous section, Theorem \ref{T:
evenDecompThm} follows.

%%% ----------------------------------------------------------------------

\section{Dual Cohomology of Symmetric Powers}

In this section, we will define a dual space $H_k^{1*}$ to the
cohomology space $H_k^1$ and a dual operator $\bar \beta_k^*$ to the
Dwork operator $\bar \beta_k$.

\bigskip\noindent{\bf Dagger Spaces.} Define
\begin{align*}
\hatR &:= \text{ analytic functions on } 1 < |a| < e, \text{
$e$ unspecified} \\
L^\dag &:= \text{ analytic functions on } |a| < e, \text{ where } e
> 1 \text{ unspecified} \\
\Ldual &:= \text{ analytic functions on } 1 <  |a| \leq \infty.
\end{align*}
Define a perfect pairing $\langle \cdot, \cdot \rangle: \hatR \times
\hatR \rightarrow \bb C_p$ by $\langle u, v \rangle := \text{ the
constant term of the product } uv$. By \cite{pR94}, this places
$\hatR$ into (topological) duality with itself. Further, by the
Mittag-Leffler theorem, $\hatR = a L^\dag \oplus \Ldual$. Thus,
since the annihilator of $L^\dag$ in $\hatR$ is $a L^\dag$, the
pairing places $L^\dag$ into duality with $\hatR / a L^\dag =
\Ldual$.

With any $b',b > 0$, define
\begin{align*}
\c M_a &:= \c M_a(b',b) \otimes_{L(b')} L^\dag \\
\c M_a^* &:= \c M_a(b',b) \otimes_{L(b')} \hatR / a L^\dag.
\end{align*}
Fixing $\{ x, x^2 \}$ as a basis of both $\c M_a$ and $\c M_a^*$, we
may place them into (topological) duality by defining the pairing
$\langle \cdot, \cdot, \rangle: \c M_a \times \c M_a^* \rightarrow
\bb C_p$ by taking $h = h_1(a) x + h_2(a) x^2 \in \c M_a$ and $g =
g_1(a) x + g_2(a) x^2 \in \c M_a^*$ and defining
\[
\langle h, g \rangle := \langle h_1, g_1 \rangle + \langle h_2, g_2
\rangle.
\]

We may extend this duality, and the pairing, to the symmetric powers
of these spaces as follows. Let $\Msk{a}$ (resp. $\MskDual{a}$) be
the $k$-th symmetric tensor power of $\c M_a$ (resp. $\c M_a^*$)
over $L^\dag$ (resp. $\Ldual$). For both $\Msk{a}$ and
$\MskDual{a}$, fix the basis $\c B := \{ v^j w^{k-j} \}_{j=0}^k$
where $v := x$ and $w := x^2$. With $(u_1, \ldots, u_k) \in (\c
M_a)^k$ and $(v_1, \ldots, v_k) \in (\c M_a^*)^k$, define a perfect
pairing between $\Msk{a}$ and $\MskDual{a}$ by taking their images
$u_1 \cdots u_k \in \Msk{a}$ and $v_1 \cdots v_k \in \MskDual{a}$
and defining
\[
( u_1 \cdots u_k, v_1 \cdots v_k) := \sum_{\sigma \in \c S_k}
\prod_{i=1}^k \langle u_i, v_{\sigma(i)} \rangle
\]
where $\c S_k$ denotes the full symmetric group on $\{1, 2, \ldots,
k\}$.

%%%-----------------------------------------------------------------------
\bigskip\noindent{\bf Differential Operator.}
With respect to the basis $\{x, x^2 \}$ of $\c M_a$, we may write
the matrix representation of $\partial_a$ as $a \frac{d}{da} - G$
with
\[
G := \left(
  \begin{array}{cc}
    0 & \pi a \\
    \frac{-\pi a^2}{3} & 0 \\
  \end{array}
\right).
\]
Note, we are thinking of $\c M_a$ as a \emph{row space}, and so, $G$
acts on the right. Define the endomorphism of $\c M_a^*$
\[
\partial^*_a := - a \frac{d}{da} - Trunc_a \circ G.
\]
The map $\partial_a^*$ acts on the \emph{column space} $\c M_a^*$.
Similar to that of $D_a$ and $D_a^*$, $\partial_a$ and
$\partial_a^*$ are dual to one another with respect to the pairing
$\langle \cdot, \cdot \rangle$.

As we did in \S \ref{S: cohomForm}, we extend $\partial_a^*$ to an
operator on $\MskDual{a}$ by extending linearly the action
\[
\partial_a^*(u_1 \cdots u_k) := \sum_{i=1}^k u_1 \cdots \hat u_i \cdots
u_k \partial_a^*(u_i)
\]
where $\hat u_i$ means we are leaving it out of the product. It is
not hard to verify that $\partial_a^*$ and $\partial_a$ are dual to
one another via the pairing $(\cdot, \cdot)$.

We may now define
\[
H_k^1 := \Msk{a} / \partial_a \Msk{a} \quad \text{and} \quad
H_k^{1*} := ker( \partial_a^* | \MskDual{a}).
\]
While we have already defined $H_k^1$ as the quotient space of $\c
M_a^{(k)}(b',b)$ by $\partial_a \c M_a^{(k)}(b',b)$, the two
definitions are the same. Since $H_k^1$ is of finite dimension, from
\cite[Proposition 7.2.2]{pR94}, $\partial_a \Msk{a}$ is a closed
subspace of $\Msk{a}$. It follows that $H_k^{1*}$ is dual to
$H_k^1$, and hence, they have the same dimension.

%%%-----------------------------------------------------------------------
\bigskip\noindent{\bf Dwork Operators.}
Similar to that of \S \ref{S: cohomForm}, define the Dwork operator
\[
\beta_k := \psi_a \circ \bar \alpha_k(a): \Msk{a} \rightarrow
\Msk{a}.
\]
Observe $\beta_k \circ \partial_a = p \partial_a \circ \beta_k$, and
so $\beta_k$ induces a map $\bar \beta_k: H_k^1 \rightarrow H_k^1$.

With respect to the ordered basis $\c B := (v^k, v^{k-1}w, \ldots,
w^k)$, we may denote by $Sym^k(\f A(a))$ the matrix of $\bar
\alpha_k(a)$. This makes the Dwork operator representable as
$\beta_k = \psi_a \circ Sym^k(\f A(a))$. From this perspective, it
is natural to define the dual of $\beta_k$ as
\[
\beta^*_k := Trunc_a \circ Sym^k(\f A(a)) \circ \Phi_a : \MskDual{a}
\rightarrow \MskDual{a}.
\]
Notice that since the matrix entries of $\f A(a)$ lie in $L^\dag$,
$\beta_k^*$ is well-defined.

Dual to $p \partial_k \circ \beta_k = \beta_k \circ
\partial_k$ is $p \beta_k^* \circ \partial_k^* = \partial_k^* \circ
\beta_k^*$. Hence, $\beta_k^*$ induces a map $\bar \beta_k^*:
H^{1*}_k \rightarrow H^{1*}_k$. It follows that
\[
det(1 - \bar \beta_k T | H^1_k) = det(1 - \bar \beta^*_k T |
H^{1*}_k).
\]

%%% ----------------------------------------------------------------------
\subsection{Primitive Cohomology and its Dual}\label{S: TrivialSubspaces}

In this section, we will decompose the (conjugate) dual
cohomology into three subspaces: a constant subspace, a trivial
subspace, and a primitive part. The eigenvectors and associated
eigenvalues of the dual Frobenius on the first two spaces are
explicitly given.

As we did in \S \ref{S: RelFunEqu}, we need to distinguish between
$\pi$ and its conjugate $-\pi$. Using notation from \S \ref{S:
RelFunEqu}, for any $b', b > 0$, define the spaces
\begin{align*}
\c M_{-\pi, a} &:= \c M_{-\pi, a}(b',b) \otimes_{L(b')} L^\dag \\
\c M_{-\pi, a}^* &:= \c M_{-\pi, a}(b',b) \otimes_{L(b')} \hatR / a
L^\dag.
\end{align*}
Denote by $\Msk{-\pi, a}$ and $\MskDual{-\pi, a}$ the $k$-th
symmetric tensor powers of these spaces over $L^\dag$ and $\Ldual$,
respectively. The differential operator on $\Msk{-\pi,a}$ is
$\partial_{-\pi, a} := a \frac{d}{da} + G_k$, and the differential
operator on $\MskDual{-\pi,a}$ is $\partial_{-\pi, a}^* := - a
\frac{d}{da} + Trunc_a \circ G_k$. From this, we may define
$H_{-\pi, k}^1 := \Msk{-\pi, a} / \partial_{-\pi, a} \Msk{-\pi,a}$
and $H_{-\pi, k}^{1*} := ker( \partial_{-\pi, a}^* | \MskDual{-\pi,
a})$.

%%%-----------------------------------------------------------------------
\bigskip\noindent{\bf Constant Subspace $\bb C_p^{k+1}$.}
Using the basis $\c B$ defined above, define $\bb C_p^{k+1} :=
\bigoplus_{j=0}^k \bb C_p v^j w^{k-j}$. Notice that $\bb C_p^{k+1}
\subset H_{-\pi, k}^{1*}$.

%%%-----------------------------------------------------------------------
\bigskip\noindent{\bf Subspace $\f T_k$ coming from infinity.} (cf.
\S \ref{S: AiryAsymp}) If $k$ is odd, define $\f T_k := 0$. Assume
$k$ is even. Define $\kappa := \frac{2i}{3\sqrt{3}}$ where $i$ is a
square root of $-1$ in $\bb C_p$, and define
\[
v_1(a) := \sum_{n=0}^\infty \frac{ \left( \frac{7}{6} \right)_n
\left( \frac{17}{12} \right)_n}{2^n \kappa^n \pi^n (n+1)!} a^{-3n},
\]
and $v_2(a)$ by replacing $i$ with $-i$ in $v_1$. By reduction,
define the vectors in $\c M_a^*$
\begin{align*}
u_1 &:= v_1 v + \left( -\frac{1}{2} a^{-2} v_1 + 3 a^2 \kappa
\pi v_1 + a^{-1} v_1' \right) w \\
u_2 &:= v_2 v + \left( -\frac{1}{2} a^{-2} v_2 - 3a^2 \kappa \pi v_2
+ a^{-1} v_2' \right) w.
\end{align*}
Next, for each $0 \leq j \leq k/2$ and $p | (k-2j)$, define
\begin{align*}
\eta_j^+ := a^{-k/2} \left( e^{(k-2j) \kappa \pi a^3} u_1^{k-j}
u_2^j +
e^{-(k-2j) \kappa \pi a^3} u_1^j u_2^{k-j}\right) \\
\eta_j^- := a^{-k/2} \left( e^{(k-2j) \kappa \pi a^3} u_1^{k-j}
u_2^j - e^{-(k-2j) \kappa \pi a^3} u_1^j u_2^{k-j}\right).
\end{align*}
If $4 | k$, then for each $0 \leq j \leq k/2$ and $p | (k-2j)$,
define $\omega_j(a)$ by $\omega_j(a^2) = \eta_j^+(a)$. Else, if $2 |
k$ but $4 \nmid k$, then for each $0 \leq j < k/2$ and $p | (k-2j)$,
define $\omega_j(a^2) = \eta_j^-(a)$.

It follows from \S \ref{S: AiryAsymp} that $\f T_k := span_{\bb
C_p}\{ Trunc_a(\omega_j(a)) \}$ is a subspace of $H_{-\pi, k}^{1*}$
of dimension
\[
dim_{\bb C_p} \f T_k =
\begin{cases}
1 + \left\lfloor \frac{k}{2p} \right\rfloor & \text{if $4 | k$} \\
\left\lfloor \frac{k}{2p} \right\rfloor & \text{if $2 | k$ but $4
\nmid k$} \\
0 & \text{if $k$ odd}.
\end{cases}
\]
To see this, from \S \ref{S: AiryAsymp} $\omega_j(a)$ is a solution
of $a \frac{d}{da} - G_k$ whose coordinates live in $\hatR$. By
Mittag-Leffler, $Trunc_a(\omega_j(a))$ is analytic on
$D^-(\infty,1)$, and while $\left(a \frac{d}{da} - G_k\right)
Trunc_a(\omega_j(a))$ is not zero, it only consists of positive
powers of the variable $a$. Hence, it is zero when we apply
$\partial_{-\pi, a}^*$.

%%%-----------------------------------------------------------------------
\bigskip\noindent{\bf Primitive Cohomology.} Define the space
\begin{equation}\label{E: subspaceDecomp}
PH_{-\pi, k}^{1*} := H_{-\pi, k}^{1*} / (\bb C_p^{k+1} \oplus \f
T_k).
\end{equation}
Using this, we define the subspace $PH_{-\pi, k}^1 \subset
H_{-\pi,k}^1$ as the annihilator of $\bb C_p^{k+1} \oplus \f T_k$ in
$H_{-\pi, k}^1$ via the pairing. Thus, $PH_{-\pi, k}^{1*}$ and
$H_{-\pi, k}^1$ are dual.

%%%-----------------------------------------------------------------------
\bigskip\noindent{\bf Frobenius on the Constant Subspace.}
Let us show
\begin{equation}\label{E: ConstAction}
\bar \beta_{-\pi,k}^*|_{\bb C_p^{k+1}} = Sym^k( \f A(0) ).
\end{equation}
To prove this, recall $\bar \beta_{-\pi,k}^* := Trunc_a \circ Sym^k
\f A_{-\pi}(a) \circ \Phi_a$. The operator $\Phi_a$ is the identity
map on $\bb C_p^{k+1}$ and $Sym^k \f A_{-\pi}(a)$ is a matrix with
entries in $L^\dag$; in particular, they are power series in the
variable $a$. Therefore, only the constant terms of these power
series will have an effect on $\bb C_p^{k+1}$, everything else would
be killed by the truncation operator. That is,
\begin{align*}
Trunc_a \circ Sym^k \f A_{-\pi}(a) \circ \Phi_a (\bb C_p^{k+1}) &=
Trunc_a \circ Sym^k \f A_{-\pi}(a) (\bb C_p^{k+1}) \\
&= Trunc_a \circ Sym^k \f A(0) (\bb C_p^{k+1})
\end{align*}
which proves (\ref{E: ConstAction}). Further, since $\f A(0)$ is
invertible, we have $\bar \beta_{-\pi, k}^*(\bb C_p^{k+1}) = \bb
C_p^{k+1}$.

%%%-----------------------------------------------------------------------
\bigskip\noindent{\bf Frobenius on the Trivial Subspace $\f T_k$.}
Recall from \S \ref{S: AiryAsymp} that
\[
\f A_{-\pi}^*(a^2) a^{-p/2} V(a^p) S(a^p) = a^{-1/2} V(a) S(a) M
\]
where $V(a)$ is the matrix $(u_1 \> \> u_2)$, $S(a) :=
\left(%
\begin{array}{cc}
  e^{\kappa \pi a^3} & 0 \\
  0 & e^{-\kappa \pi a^3} \\
\end{array}%
\right)$, and $M$ is a constant matrix given by
\[
M=
\left(%
\begin{array}{cc}
  \sqrt{p} & 0 \\
  0 & \pm \sqrt{p} \\
\end{array}%
\right) \quad \text{if } p \equiv 1,7 \text{ mod}(12)
\]
(of course, $\sqrt{p}$ may not be the positive square root) and
\[
M =
\left(%
\begin{array}{cc}
  0 & \sqrt{-g} \\
  \pm \sqrt{-g} & 0 \\
\end{array}%
\right) \quad \text{if } p \equiv 5,11 \text{ mod}(12)
\]
where ``$\pm$'' means positive if $p \equiv 1, 5$ mod(12) and
negative if $p \equiv 7,11$ mod(12), and $g := g_2((p^2-1)/3)$.

From this, we may calculate that for $4 | k$
\[
\bar \beta_{-\pi, k}^*( Trunc_a(\omega_j) ) =
\begin{cases}
p^{k/2} Trunc_a(\omega_j) & \text{if } p \equiv 1 \text{mod(12)} \\
(-1)^j p^{k/2} Trunc_a(\omega_j) & \text{if } p \equiv 7
\text{mod(12)} \\
g^{k/2} Trunc_a(\omega_j) & \text{if } p \equiv 5 \text{
mod(12)}
\\
(-1)^j g^{k/2} Trunc_a(\omega_j) & \text{if } p \equiv 11
\text{ mod(12)}.
\end{cases}
\]
and for $2|k$ but $4 \nmid k$,
\[
\bar \beta_{-\pi, k}^*( Trunc_a(\omega_j) ) =
\begin{cases}
p^{k/2} Trunc_a(\omega_j) & \text{if } p \equiv 1 \text{mod(12)} \\
(-1)^j p^{k/2} Trunc_a(\omega_j) & \text{if } p \equiv 7
\text{mod(12)} \\
(-g)^{k/2} Trunc_a(\omega_j) & \text{if } p \equiv 5 \text{
mod(12)} \\
(-1)^{j+1} (-g)^{k/2} Trunc_a(\omega_j) & \text{if } p
\equiv 11 \text{ mod(12)}.
\end{cases}
\]
This demonstrates that not only is $\bar \beta_{-\pi, k}^*$
stable on $\f T_k$, but the basis $\{ Trunc_a(\omega_j) \}$
consists of eigenvectors.

Using either Gauss sums or Dwork theory, one may show
\[
P_{-\pi, k}(T) := det(1-\bar \beta_{-\pi,k}^* T | \bb C_p^{k+1}) =
\prod_{j=0}^k (1 - \pi_1(0)^j \pi_2(0)^{k-j} T)
\]
where $L(x^3 / \bb F_p, T) = (1 - \pi_1(0) T)(1 - \pi_2(0) T)$;
note, we are using the (conjugate) splitting function $\exp(-\pi(t -
t^p))$ to compute this $L$-function. Letting $N_{-\pi,k}(T) :=
det(1-\bar \beta_{-\pi,k}^* T| \f T_k)$, we see that
\[
\overline{M}_k^*(T) = \frac{det(1-\bar \beta_{-\pi,k}^* T |
H_{-\pi,k}^{1*})}{det(1 - p \bar \beta_{-\pi, k} T | H_{-\pi,k}^0)}
= \frac{P_{-\pi,k}(T) N_{-\pi,k}(T) det(1-\bar \beta_{-\pi,k} T |
PH_{-\pi,k}^1)} {det(1 - p \bar \beta_{-\pi, k} T | H_{-\pi,k}^0)}
\]
where we have used duality for the second equality and the bar on
the left means complex conjugation. Note, conjugation acts as the
identity on $M_k^*(T)$ since it has real coefficients. It follows
from (\ref{E: Relations}) that
\[
M_k(T) = \frac{N_{\pi,k}(T) det(1-\bar \beta_{\pi,k} T |
PH_{\pi,k}^1)} {det(1 - p \bar \beta_{\pi, k} T | H_{\pi,k}^0)}
\]

%%% ----------------------------------------------------------------------
\subsection{Functional Equation}
In this section, we will define an isomorphism $\Theta_k : PH_{-\pi,
k}^{1*} \rightarrow PH_{\pi, k}^1$ which relates the Frobenius $\bar
\beta_{\pi, k}$ with its conjugate dual $\bar z_{-\pi, k}^*$; see
equation (\ref{E: realFunEqu}).

For convenience, denote by $\f A_{-\pi, k, a}$ the matrix $Sym^k(\f
A_{-\pi}(a))$. Dual to $\beta_{-\pi,k} := \psi_a \circ \f A_{-\pi,
k, a}: \Msk{-\pi,a} \rightarrow \Msk{-\pi,a}$ is the isomorphism
$z_{-\pi,k}^* := Trunc_a \circ \f A_{-\pi,k,a} \circ \Phi_a:
\MskDual{-\pi,a} \rightarrow \MskDual{-\pi,a}$. For $\xi^* \in
ker(\partial_{-\pi,a}^* | \MskDual{-\pi,a})$, we have (in $\bb
C_p[[a^{\pm 1}]]^{k+1}$)
\begin{equation}\label{E: FunBegin}
\f A_{-\pi, k,a} \circ \Phi_a(\xi^*) = z_{-\pi,k}^*(\xi^*) + \eta
\end{equation}
for some $\eta \in \bb C_p[a]^{k+1}$.

Recall the space $\c R_{-\pi, a}'(b',b)$ from \S \ref{S: RelFunEqu}.
Define $\c R_{-\pi, a}' := \c R_{-\pi, a}'(b',b) \otimes_{L(b')}
L^\dag$. Define
\[
\Theta_{k,a}' := -a \frac{d}{da} + G_k: ker(\partial_{-\pi, a}^* |
\MskDual{-\pi, a}) \rightarrow Sym^k( \c R_{-\pi, a}').
\]
It will be useful to note that $\Theta_{k,a}'$ is both
$\partial_{-\pi,a}^*$ without the truncation map (and hence it
$p$-commutes with $\f A_{-\pi,k,a} \circ \Phi_a$) and
$-\partial_{\pi,a}$. Applying $\Theta_{k,a}'$ to (\ref{E:
FunBegin}), we have (in $Sym^k(R_{-\pi,a}')$)
\begin{equation}\label{E: Fun2}
p \> \f A_{-\pi,k,a} \circ \Phi_a \circ \Theta_{k,a}'(\xi^*) =
\Theta_{k,a}' \circ z_{-\pi, k}^*(\xi^*) + \Theta_{k,a}'(\eta).
\end{equation}

Recall from \S \ref{S: RelFunEqu} the isomorphism
$\overline{\Theta}_{-\pi, a}: \c R_{-\pi,a}'(b',b^*) \rightarrow \c
M_{\pi, a}(b',b)$ which satisfies $\overline{\Theta}_{-\pi, a^p} =
p^{-1} \bar \alpha_\pi(a) \circ \overline{\Theta}_{-\pi, a} \circ
\bar \alpha_{-\pi}^*(a)$. Denoting by $\overline{\Theta}_{-\pi,k,a}$
the $k$-th symmetric power of $\overline{\Theta}_{-\pi, a}$, and
using that the matrix of $\bar \alpha_{-\pi}^*(a)$ equals $\f
A_{-\pi}(a)$, we have
\begin{equation}\label{E: SymRelFunEqu}
p \f A_{\pi, k, a}^{-1} \circ \overline{\Theta}_{-\pi, k, a^p} =
\overline{\Theta}_{-\pi,k, a} \circ \f A_{-\pi, k, a}
\end{equation}
Thus, applying $\overline{\Theta}_{-\pi, k, a}$ to both sides of
(\ref{E: Fun2}), we have
\[
p^{k+1} \f A_{\pi, k, a}^{-1} \circ \overline{\Theta}_{-\pi, k, a^p}
\circ \Phi_a \circ \Theta_{k,a}'(\xi^*) = \overline{\Theta}_{-\pi,
k, a} \Theta_{k,a}' z_{-\pi, k}^*(\xi^*) + \overline{\Theta}_{-\pi,
k, a} \Theta_{k,a}'(\eta)
\]
which equals
\[
p^{k+1} \f A_{\pi, k,a}^{-1} \circ \Phi_a \circ
\overline{\Theta}_{-\pi, k, a} \Theta_{k,a}'(\xi^*) =
\overline{\Theta}_{-\pi, k, a} \Theta_{k,a}' z_{-\pi, k}^*(\xi^*) +
\overline{\Theta}_{-\pi, k, a} \Theta_{k,a}'(\eta).
\]
Thus,
\begin{equation}\label{E: almostfunequ}
\overline{\Theta}_{-\pi, k, a} \Theta_{k,a}'(\xi^*) = p^{-(k+1)}
\beta_{\pi,k} \overline{\Theta}_{-\pi, k, a} \Theta_{k,a}' z_{-\pi,
k}^*(\xi^*) + \beta_{\pi,k} \overline{\Theta}_{-\pi, k, a}
\Theta_{k,a}'(\eta).
\end{equation}
Now,
\begin{align*}
\beta_{\pi,k} \overline{\Theta}_{-\pi, k, a} \Theta_{k,a}'(\eta) &=
\psi_a \circ p \overline{\Theta}_{-\pi, k, a^p} (\bar \alpha_{-\pi,
k, a}^*)^{-1} \Theta_{k,a}'( \eta) \\
&= p \overline{\Theta}_{-\pi, k, a} \circ \psi_a \circ \Theta_{k,
a^p}' (\bar \alpha_{-\pi, k, a}^*)^{-1}(\eta) \\
&= p^2 \Theta_{k, a}' \overline{\Theta}_{-\pi, k, a} \psi_a \circ
(\bar \alpha_{-\pi, k, a}^*)^{-1}(\eta) \\
&\equiv 0 \text{ mod}(\partial_{\pi,a} \Msk{\pi,a})
\end{align*}
where the first equality comes from (\ref{E: SymRelFunEqu}).
Finally, define $\Theta_k : ker(\partial_{-\pi, k}^* |
\MskDual{-\pi,a}) \rightarrow H_{\pi,k}^1$ by
\[
\Theta_k(\xi^*) := \overline{\Theta}_{-\pi, k, a}
\Theta_{k,a}'(\xi^*) \quad \text{ mod}(\partial_{\pi,a}
\Msk{\pi,a}).
\]
Then (\ref{E: almostfunequ}) shows us that
\[
p^{k+1} \bar \beta_{\pi, k}^{-1} \Theta_k = \Theta_k \bar z_{-\pi,
k}^*.
\]
However, this is not quite the functional equation since the map
$\Theta_k$ has a kernel.

\begin{lemma}
$ker \> \Theta_k = \bb C_p^{k+1} \oplus \f T_k.$
\end{lemma}

\begin{proof}
Since $z_{-\pi, k}^*(\bb C_p^{k+1}) = 0$, we see by (\ref{E:
almostfunequ}) that $\Theta_k(\bb C_p^{k+1}) = 0$. Next, since
$\Theta_{k,a}'$ is simply $\partial_{-\pi, k}^*$ without the
truncation map, we have $\Theta_{k,a}'( \f T_k) = 0$, and so,
$\Theta_k(\f T_k) = 0$.

Suppose $\xi^* \in \MskDual{-\pi,a}$ is in the kernel of $\Theta_k$
and is not constant. Then there exists $h \in \Msk{\pi,a}$ such that
$\Theta_k(\xi^*) = \partial_{\pi, a}(h)$. Now,
$\overline{\Theta}_{-\pi, k, a} \Theta_{k,a}' = \Theta_{k,a}'
\overline{\Theta}_{-\pi,k,a}$, and so
\[
\left( a \frac{d}{da} - G_k \right)(-\xi^* + \tilde h) = 0
\]
where $\tilde h := \overline{\Theta}_{-\pi,k,a}^{-1}(h)$. Thus,
$-\xi^* + \tilde h$ corresponds to a solution in $ker(l_k | \hatR)$.
But these solutions are in one-to-one correspondence with $\f T_k$.
\end{proof}

Thus, if we restrict to $\bar z_{-\pi, k}^*$ on $PH_{-\pi, k}^{1*}$
and $\bar \beta_{\pi, k}$ on $PH_{\pi, k}^1$, then $\Theta_k:
PH_{-\pi, k}^{*1} \rightarrow PH_{\pi, k}^1$ is an isomorphism which
satisfies
\begin{equation}\label{E: realFunEqu}
p^{k+1} \bar \beta_{\pi, k}^{-1} \Theta_k = \Theta_k \bar z_{-\pi,
k}^*.
\end{equation}

\begin{corollary}
$dim_{\bb C_p} PH_k^{1} = dim_{\bb C_p} PH_{-\pi,k}^{1*} \leq k$.
\end{corollary}

\begin{proof}
From the explicit description of the matrix $G_k$, the image of the
\emph{injective} function $\Theta_k': PH_{-\pi,k}^{1*} \rightarrow
\bb C_p[a]^{k+1}$ is contained within the subspace $\{0 \} \times a
\bb C_p^k$ which has dimension $k$.
\end{proof}

%%% ----------------------------------------------------------------------
\section{Newton Polygon of $M_k(T)$}\label{S: NPofM}

For $k$ odd and $k < p$, from Theorem \ref{T: oddDecomp} and
Section \ref{S: TrivialSubspaces}, we have
\[
M_k(T) = det(1-\bar \beta_k T| PH_k^1)
\]
where
\[
PH_k^1 = \bigoplus_{j=0}^{(k-1)/2} (\bb C_p a) v^{k-2j} w^{2j}.
\]
Using this explicit description of the cohomology, and the
effective decomposition theorem (Theorem \ref{T: oddDecomp}),
we may prove the following.

\begin{theorem}\label{T: OddNewtonPolygon}
Suppose $k$ is odd and $k < p$. Writing
\[
M_k(T) = 1 + c_1 T + c_2 T^2 + \cdots
+c_{(k+1)/2} T^{(k+1)/2}
\]
we have
\[
ord(c_m) \geq \frac{(p-1)^2}{3p^2}(m^2 + m + mk).
\]
for every $m = 0, 1, \ldots, (k+1)/2$.
\end{theorem}

\begin{proof}
Let $b := (p-1)/p$ and $b' := b/p$. From \S \ref{S: FrobEstimates},
we know $\alpha(a) x^i \equiv \f A_{i1}(a) x + \f A_{i2}(a) x^2$
mod($D_a \c K(b', b)$) with $\f A_{ij} \in L(b';
\frac{b'}{3}(pj-i))$. From \S \ref{S: oddCohom}, we know $\{ a
v^{k-2i} w^{2i} \}_{i=0}^{(k-1)/2}$ is a basis of $PH_k^1$. Since $v
:= x$ and $w := x^2$,
\[
\bar \alpha_k(a)( a v^{k-2i} w^{2i} ) = a (\bar \alpha(a)v)^{k-2i}
(\bar \alpha(a)w)^{2i} = \sum_{r=0}^{k} a B_r(a) v^{k-r} w^{r}
\]
where
\[
B_{r}(a) := \sum_{\substack{n + m = k-r \\ n=0,1, \ldots, k-2i \\
m=0,1,\ldots, 2i}} \binom{k-2i}{n} \binom{2i}{m} (\f A_{11})^n (\f
A_{21})^{k-2i-n} (\f A_{12})^{m} (\f A_{22})^{2i-m}.
\]
After some calculation, we see that
\[ a B_{2r} \in L(b';
\frac{b'}{3}[(p-1)k + pr - 2i]- 2b'/3) = L(b'; \frac{b'}{3}[(p-1)k -
2i] + \delta_{r+1}-\frac{b}{3} - \frac{2b'}{3}).
\]
where $\delta_{r+1} := \frac{b(r+1)}{3}$. Using notation from \S
\ref{S: oddCohom}, this means
\[
\bar \alpha_k(a)( a v^{k-2i} w^{2i}) \in \Msk{a^p}(b';\bs \delta +
\frac{b'}{3}[(p-1)k - 2i] - \frac{b}{3} - \frac{2b'}{3}).
\]
Consequently,
\[
\beta_k( a v^{k-2i} w^{2i}) \in \Msk{a}(b; \bs \delta +
\frac{b'}{3}[(p-1)k - 2i] - \frac{b}{3} - \frac{2b'}{3}).
\]
It follows from Lemma \ref{L: SymDecomp_D} that for some
constants $A_{ij}$ we may write
\[
\beta_k( a v^{k-2i} w^{2i}) \subset \sum_{j=0}^{(k-1)/2} A_{ij} a
v^{k-2j} w^{2j} + \partial_a \Msk{a}(b)
\]
with the sum being an element of $V_k(b;\bs \delta +
\frac{b'}{3}[(p-1)k - 2i] - \frac{b}{3} - \frac{2b'}{3})$.
Therefore,
\[
A_{ij}a \in L(b; \delta_{2j+1} + \frac{b'}{3}[(p-1)k - 2i] -
\frac{b}{3} - \frac{2b'}{3}),
\]
and so,
\[
ord(A_{ij}) \geq \frac{2b'}{3}(pj-i) + \frac{b'}{3}(p-1)k +
\frac{2}{3}(b-b').
\]
Let us rewrite this. Let $\xi \in \bb C_p$ such that $ord_p(\xi) =
2b'/3$. Fix $\tilde \xi := \xi^{p-1}$. Then, if we had used the
basis $\{ \xi^i a v^{k-2i} w^{2i} \}_{i=0}^{(k-1)/2}$, then
\[
\bar \beta_k( \xi^i a v^{k-2i} w^{2i} ) = \sum_{j=0}^{(k-1)/2}
\left(\xi^{i-j} A_{ij} \right) \xi^j a v^{k-2j} w^{2j}.
\]
Hence, the matrix of $\bar \beta_k$ on $PH_k^1$, with respect to
this new basis, takes the form
\begin{equation}\label{E: Filtration}
\text{matrix of } \bar \beta_k = p^{\frac{b'}{3}(p-1)k +
\frac{2}{3}(b-b')} \left(
  \begin{array}{cccc}
    B_{00} & \tilde \xi B_{01} & \cdots & \tilde\xi^{(k-1)/2} B_{0, (k-1)/2} \\
    B_{10} & \tilde \xi B_{11} & \cdots & \tilde\xi^{(k-1)/2} B_{1, (k-1)/2} \\
    \vdots & \vdots & \vdots & \vdots \\
    B_{(k-1)/2, 0} & \tilde \xi B_{(k-1)/2, 1} & \cdots & \tilde\xi^{(k-1)/2} B_{(k-1)/2, (k-1)/2} \\
  \end{array}
\right),
\end{equation}
where $\tilde \xi^j B_{ij} = \xi^{i-j} A_{ij}$, and so,
$ord_p(B_{ij}) \geq 0$. It follows that, if we write
\[
det(1- \bar \beta_k T| PH_k^1) = 1 + c_1 T + \cdots +c_{(k+1)/2}
T^{(k+1)/2}
\]
then since $b := (p-1)/p$ and $b' := b/p$,
\[
ord(c_m) \geq \frac{(p-1)^2}{3p^2} ( m^2+ m+ mk ).
\]
\end{proof}

%%%-----------------------------------------------------------------------
\bibliographystyle{amsplain}
\bibliography{XBib}

\providecommand{\bysame}{\leavevmode\hbox to3em{\hrulefill}\thinspace}
\providecommand{\MR}{\relax\ifhmode\unskip\space\fi MR }
\providecommand{\MRhref}[2]{%
  \href{http://www.ams.org/mathscinet-getitem?mr=#1}{#2}
}
\providecommand{\href}[2]{#2}
\begin{thebibliography}{10}

\bibitem{AdolphHecke}
A.~Adolphson, \emph{A $p$-adic theory of {H}ecke polynomials}, Duke Math.
  Journal \textbf{43} (1976), no.~1, 115--145.

\bibitem{BlacheFerard}
R.~Blache and \'E. F\'erard, \emph{{N}ewton stratification for polynomials: The
  open stratum}, J. Number Th. \textbf{123} (2007), no.~2, 456--472.

\bibitem{ZhuBlacheFerard}
R.~Blache, \'E. F\'erard, and J.~Zhu, \emph{{H}odge-{S}tickelberger polygons
  for {$L$}-functions of exponential sums of $p(x^s)$}.

\bibitem{bD60}
B.~Dwork, \emph{On the rationality of the zeta function of an algebraic
  variety}, Amer. J. of Math. \textbf{82} (1960), no.~3, 631 -- 648.

\bibitem{DwZetaIHES}
\bysame, \emph{On the zeta function of a hypersurface}, Pub. Math. I.H.E.S.,
  Paris (1962), no.~12, 5 -- 68.

\bibitem{DwIII}
\bysame, \emph{On the zeta function of a hypersurface: {III}}, Annals of Math
  \textbf{83} (1966), no.~3, 457--519.

\bibitem{bD69}
\bysame, \emph{$p$-adic {C}ycles}, Pub. Math. I.H.E.S., Paris \textbf{37}
  (1969), p. 27--115.

\bibitem{DwHecke}
\bysame, \emph{On {H}ecke polynomials}, Inventiones math. \textbf{12} (1971),
  249--256.

\bibitem{bD80}
\bysame, \emph{{L}ectures on $p$-adic differential equations}, Springer-Verlag,
  1980.

\bibitem{FuWan}
L.~Fu and D.~Wan, \emph{{$L$}-functions for symmetric products of {K}loosterman
  sums}, J. Reine Angew. Math. \textbf{589} (2005), 79 -- 103.

\bibitem{LivneCubic}
R.~Livn\'{e}, \emph{The average distribution of cubic exponential sums}, J.
  Reine Angew. Math. \textbf{375/376} (1987), 362--379.

\bibitem{pM70}
P.~Monsky, \emph{$p$-adic analysis and zeta functions}, Lectures in Math.,
  Kinokuniya Book-Store, Japan, 1970.

\bibitem{pR86}
P.~Robba, \emph{Symmetric powers of the $p$-adic {B}essel equation}, J. Reine
  Angew. Math. \textbf{366} (1986), 194 -- 220.

\bibitem{pR94}
P.~Robba and G.~Christol, \emph{{\'E}quations diff\'erentielles $p$-adiques},
  Hermann, 1994, French.

\bibitem{SperberCubic}
S.~Sperber, \emph{On the $p$-adic theory of exponential sums}, Amer. J. Math.
  \textbf{108} (1986), 255--296.

\bibitem{WanModular}
D.~Wan, \emph{Dimension variation of classical and $p$-adic modular forms},
  Invent. Math. \textbf{133} (1998), 469--498.

\bibitem{WanVariationNP}
\bysame, \emph{{V}ariation of $p$-adic {N}ewton polygons for {$L$}-functions of
  exponential sums}, Asian J. Math. \textbf{8} (2004), no.~3, 427--474.

\bibitem{Yang}
R.~Yang, \emph{{N}ewton polygons of {$L$}-functions of polynomials of the form
  $x^d+\lambda x$}, Finite Fields and App. \textbf{9} (2003), 59--88.

\bibitem{ZhuAiry}
J.~Zhu, \emph{$p$-adic variation of {$L$}-functions of one variable exponential
  sums, i}, Amer. J. Math. \textbf{125} (2003), 669 -- 690.

\bibitem{ZhuAsympVariation}
\bysame, \emph{Asymptotic variation of {$L$}-functions of one-variable
  exponential sums}, J. Reine Angew. Math. \textbf{572} (2004), 219--233.

\end{thebibliography}

\bigskip\noindent
{\it C. Douglas Haessig \\
Department of Mathematics \\
University of Rochester \\
Rochester, NY 14627}

\end{document}